\documentclass[a4paper, 11pt]{amsart}
\usepackage[left=30truemm,right=30truemm,top=30truemm,bottom=30truemm]{geometry}

\usepackage{float}
\usepackage[colorlinks=True, citecolor=blue, linkcolor=blue]{hyperref}
\usepackage{comment}
\usepackage{amssymb}
\usepackage{tikz-cd}
\usepackage{tkz-graph}
\usepackage[capitalize, nameinlink]{cleveref}
\usepackage{booktabs}
\usepackage{multirow}
\usepackage{mathrsfs} 
\usepackage{mathtools} 
\usepackage{bm}
\usepackage{float} 
\usepackage[normalem]{ulem} 
\usepackage{xcolor}
\usepackage[all]{xy} 
\xyoption{poly}
\usepackage {units} 

\theoremstyle{definition}
\newtheorem{definition}{Definition}[section]
\newtheorem{example}[definition]{Example}
\newtheorem{remark}[definition]{Remark}

\theoremstyle{plain}
\newtheorem{theorem}[definition]{Theorem}
\newtheorem{lemma}[definition]{Lemma}
\newtheorem{proposition}[definition]{Proposition}
\newtheorem{corollary}[definition]{Corollary}
\newtheorem{conjecture}[definition]{Conjecture}
\newtheorem{algorithm}[definition]{Algorithm}

\newenvironment{steps}{%
  \begin{list}{\textup{Step \arabic{enumi}.}}{%
    \usecounter{enumi}%
    \settowidth{\labelwidth}{\textup{Step 10.}}%
    \setlength{\leftmargin}{\labelwidth}%
    \addtolength{\leftmargin}{\labelsep}%
  }%
}{\end{list}}

\DeclareMathOperator{\Tr}{Tr}

\DeclareMathOperator{\rk}{rk}

\DeclareMathOperator{\id}{id}

\DeclareMathOperator{\Hom}{Hom}
\DeclareMathOperator{\NS}{NS}

\DeclareMathOperator{\scale}{\mathfrak{s}}

\DeclareMathOperator{\im}{im}

\DeclareMathOperator{\Pic}{Pic}
\DeclareMathOperator{\Gal}{Gal}
\DeclareMathOperator{\Spec}{Spec}
\DeclareMathOperator{\Nef}{Nef}
\newcommand{\rO}{\mathrm{O}}
\newcommand{\cI}{\mathcal{I}}

\newcommand{\FF}{\mathbb{F}}

\newcommand{\ZZ}{\mathbb{Z}}

\renewcommand{\O}{\mathcal{O}}

\newcommand{\RN}[1]{\textup{\uppercase\expandafter{\romannumeral#1}}}

\usepackage{etoolbox} 
\newcommand{\bAlphabet}{
  \def\do##1{\expandafter\newcommand\csname b##1\endcsname{\mathbb{##1}}}
  \docsvlist{C,F,P,Q,R,Z}
}
\bAlphabet
\newcommand{\cAlphabet}{
  \def\do##1{\expandafter\newcommand\csname c##1\endcsname{\mathcal{##1}}}
  \docsvlist{A,C,D,E,F,K,L,O,S,T,W, G, g}
}
\cAlphabet
\newcommand{\frakAlphabet}{
  \def\do##1{\expandafter\newcommand\csname frak##1\endcsname{\mathfrak{##1}}}
  \docsvlist{a,b,D,F,f,g,m,p,P,G}
}
\frakAlphabet

\makeatletter
\newcommand{\bigperp}{%
  \mathop{\mathpalette\bigp@rp\relax}%
  \displaylimits
}

\newcommand{\bigp@rp}[2]{%
  \vcenter{
    \m@th\hbox{\scalebox{\ifx#1\displaystyle2.1\else1.5\fi}{$#1\perp$}}
  }%
}
\makeatother

\renewcommand{\o}{{\bZ[\hat \alpha]}}
\renewcommand{\O}{\mathcal{O}}

\makeatletter
\newcommand*{\defeq }{\mathrel{\rlap{%
                     \raisebox{0.3ex}{$\m@th\cdot$}}%
                     \raisebox{-0.3ex}{$\m@th\cdot$}}%
                     =}
\makeatother

\newcommand{\II}{\mathrm{II}}

\usepackage{xspace}
\newcommand{\condAst}{\ensuremath{\hyperlink{ast}{(\ast)}}\xspace}

\makeatletter
\def\blfootnote{\xdef\@thefnmark{}\@footnotetext}

\title[]{Pointless K3 surfaces over finite fields}

\author{Simon Brandhorst}
\address{Simon Brandhorst,
Fakult\"at f\"ur Mathematik und Informatik, Universit\"at des Saarlandes, Campus E2.4, 66123 Saarbr\"ucken, Germany}
\email{brandhorst@math.uni-sb.de}

\author{Yuta Takada}
\address{Yuta Takada,
Mathematical Sciences, The University of Tokyo, 
Tokyo 153-8914, Japan; JSPS Research Fellow
\vspace{1mm}
\newline
\hspace*{4mm}
Fakult\"at f\"ur Mathematik und Informatik, Universit\"at des Saarlandes, Campus E2.4, 66123 Saarbr\"ucken, Germany}
\email{ytakada@ms.u-tokyo.ac.jp}

\begin{document}

\begin{abstract}
We provide an algorithm for enumerating ordinary K3 surfaces over a fixed finite field
based on Taelman's theorem.  
We use this algorithm to prove that there exist no K3 surfaces without rational points 
over $\bF_q$ for every prime power $q\geq 13$. 
For elliptic K3 surfaces, we give a geometric proof of this 
non-existence without computer assistance. 
On the other hand, we also provide explicit pointless K3 surfaces over $\bF_q$ 
for every prime power $q\leq 11$. 
\end{abstract}

\maketitle


\section{Introduction}
Let $X$ be a K3 surface over a finite field $k$ of characteristic $p$. Assume that $X$ is \emph{ordinary}, which is equivalent to saying $\#X(k) \not\equiv 1 \bmod p$. Using the canonical lift $X_\mathrm{can}$ and an embedding $\iota: W(k)\to \bC$, where $W(k)$ is the ring of Witt vectors, one can associate a complex K3 surface $X_{\mathrm{can}}^\iota \coloneq X_{\mathrm{can}}\otimes_{W(k)} \bC$.
Moreover, the geometric Frobenius of $X$ induces a similitude $F$ of the lattice 
$\Lambda \coloneq H^2(X_{\mathrm{can}}^\iota, \bZ)$ with multiplier $q^2$ where $q = |k|$, 
and $F$ preserves the ample cone $\cK\subset \Lambda\otimes \bR$.
In his paper \cite{Taelman20}, under a mild assumption on potential semi-stable reduction, 
L. Taelman shows that the correspondence $X\mapsto (\Lambda, F, \cK)$ gives rise to an equivalence 
between the groupoid of ordinary K3 surfaces over $k$ and a groupoid of triples consisting of 
a K3 lattice $\Lambda$, an endomorphism $F$ of $\Lambda$, and a cone 
$\cK \subset \Lambda\otimes \bR$ satisfying suitable conditions.
This paper provides an algorithm for enumerating the triples arising from ordinary K3 surfaces over $k$.

For such a triple $(\Lambda, F, \cK)$, the N\'eron--Severi lattice of the corresponding K3 surface 
is recovered as 
\[ \NS(\Lambda,F) 
\coloneq \{ x\in \Lambda \mid F^dx = q^dx \text{ for some $d\geq 1$} \}. 
\]
Note that $\NS(\Lambda,F)$ is invariant under $F$, and hence, the 
transcendental lattice $T(\Lambda, F) \coloneq \NS(\Lambda,F)^\perp$ is 
also invariant under $F$. 
In practice, we will often handle $q^{-1}F:\Lambda_\bQ\to \Lambda_\bQ$ 
instead of $F$ itself since it is an isometry, although it need not preserve the lattice $\Lambda$.

The outline of the algorithm is as follows; see \S\ref{sec:strategy} for the details. 
We first fix polynomials $S$ and $C$ that will be the characteristic polynomials 
of $q^{-1}F|_{T(\Lambda, F)}$ and $q^{-1}F|_{\NS(\Lambda, F)}$, respectively. 
We then compile a list of all possible candidates for the pairs
$(T(\Lambda, F), F|_{T(\Lambda, F)})$ and $(\NS(\Lambda, F), F|_{\NS(\Lambda, F)})$. 
Finally, for each candidate for the transcendental part and the  N\'eron--Severi part, 
we glue them together to obtain a K3 lattice $\Lambda$ with similitude $F$ such that 
$F$ preserves a cone $\cK \subset \Lambda\otimes \bR$ satisfying suitable conditions.
The idea of gluing is similar to that of gluing for lattices with isometry, which can be 
found for example in \cite{McMullen2011}, but we need to be careful 
since we deal with similitudes.  

We use this algorithm to rule out pointless K3 surfaces over certain finite fields.
Here, a variety $X$ over a field $k$ is said to be \emph{pointless} if $X$ has no $k$-point.
We recall a conjecture concerning pointless varieties over finite fields.
In \cite{Bernasconi-Stefano2025}, generalizing a conjecture of
Z. Patakfalvi and M. Zdanowicz \cite{Patakfalvi-Zdanowicz2019},
F. Bernasconi and S. Filipazzi proposed the following conjecture.

\begin{conjecture}\label{conj:BS}
For every positive integer $d$, there exists an integer $q_d$ such that 
the following holds: if $X$ is a geometrically integral smooth projective
variety of dimension $d$ over $\bF_q$ and $-K_X$ is nef, and if 
$q \geq q_d$, then $X(\bF_q) \neq \emptyset$.
\end{conjecture}

We can take $q_1 = 2$, which follows, for example, from the Weil conjectures.
Let us now consider the case $d = 2$.
Let $X$ be a geometrically integral smooth projective surface 
over a finite field $k$ with $-K_X$ nef.
If $X$ is not a K3 surface, then $X(k)\neq \emptyset$;
see \cite[\S 3]{Bernasconi-Stefano2025}.
In other words, K3 surfaces are the only possible obstruction to taking $q_2 = 2$.
The Fermat quartic surface in $\bP^3$ over $\bF_5$ is a pointless K3 surface. 
Hence, we have $q_2 > 5$. On the other hand, 
by the Weil conjectures, one can take $q_2 = 23$, although this bound may not be optimal.
In dimension $3$, Bernasconi and Filipazzi show that one can take $q_3 = 23$
if the range of varieties $X$ in \cref{conj:BS} is restricted;
see the introduction of \cite{Bernasconi-Stefano2025}.
Here, the value $23$ arises from the known bound $q_2=23$.

These results motivate the following question: for which $q$ does there exist a pointless K3 surface over $\bF_q$?

\begin{theorem}\label{th:has-rational-pts}
Let $q$ be a power of a prime. If $q\geq 13$, then every K3 surface over $\bF_q$ has an $\bF_q$-point.
\end{theorem}

By contrast, for each prime power $q\leq 11$, we provide an explicit projective
model of a pointless K3 surface. As a result, we have:

\begin{theorem}\label{th:pointlessK3surfaces}
Let $q$ be a power of a prime. If $q\leq 11$, then there exists a pointless K3 surface over $\bF_q$.
\end{theorem}

These two theorems imply that the optimal prime power for $q_2$ in \cref{conj:BS} 
is $13$. 

Part of the proof of \cref{th:has-rational-pts} relies on computer calculations.
However, the following sufficient criteria for the existence of
rational points can be established without computer assistance.

\begin{theorem}\label{th:criteria-for-having-pts}
Let $X$ be a K3 surface over a finite field $k$. 
\begin{enumerate}
\item If there exists an element $\delta\in \Pic(X)$ with $\delta^2 = -2$, 
then $X(k)\neq \emptyset$. 
\item If $|k|\geq 13$ and there exists a nonzero element $\delta\in \Pic(X)$ 
with $\delta^2 = 0$, then $X(k)\neq \emptyset$. 
\end{enumerate}
\end{theorem}

The proof of \cref{th:criteria-for-having-pts} uses the geometry of 
K3 surfaces over a non-closed field. For example, a K3 surface $X$ admits a
genus-one fibration if and only if $\Pic(X)$ contains a nonzero element
$\delta$ with $\delta^2=0$ (see \cref{prop:genus-one-fibration}).
Hence, the second assertion of \cref{th:criteria-for-having-pts} shows that, 
for K3 surfaces admitting a genus-one fibration, 
\cref{th:has-rational-pts} can be proved without computer assistance.

\begin{remark}
Pointless K3 surfaces were independently considered by Auel, Engel, Kedlaya and Shaikh. They obtained  \Cref{th:pointlessK3surfaces,th:has-rational-pts} with different methodology. The AI-assisted proofs of \cref{th:criteria-for-having-pts} are similar. See their introduction \cite{AEKS} for a comparison of the two approaches.
\end{remark}

The organization of this paper is as follows. 
In \S\ref{sec:Taelman's-theorem}, we recall Taelman's theorem.
In \S\ref{sec:Frob-lattices}, we introduce \emph{Frobenius lattices} and study their properties.
A Frobenius lattice consists of a lattice and its similitude satisfying some conditions, and 
the K3 lattice equipped with the similitude $F$ described above is an example. 
As an irreducible factor of the characteristic polynomial of the isometry associated with 
such a similitude, a \emph{$q$-fractional Salem polynomial} appears. 
We study $q$-fractional Salem polynomials in \S\ref{sec:qSalem} and give 
an algorithm for enumerating them with prescribed degree. 
In \S\ref{sec:strategy}, we describe the algorithm for
enumerating the triples $(\Lambda, F, \cK)$ arising from ordinary K3 surfaces via
Taelman's theorem. The details for 
the transcendental part and the N\'eron--Severi part of $(\Lambda, F, \cK)$ are 
treated in \S\ref{sec:transcendental-part} and \S\ref{sec:Neron--Severi-part}, respectively. 
The proof of \cref{th:has-rational-pts} is also given in \S\ref{sec:Neron--Severi-part}.
In \S\ref{sec:geometric-approach}, we prove \cref{th:criteria-for-having-pts}.
Finally, we prove \cref{th:pointlessK3surfaces} in \S\ref{sec:examples} by 
giving explicit examples. 

\subsection*{Use of AI}
Generative AI assisted the authors in exploring a proof of \cref{th:criteria-for-having-pts} and in identifying some of the examples of 
pointless K3 surfaces. It was moreover used for coding assistance and for correcting the grammar and polishing the style of the manuscript. 
The authors take full responsibility for its content.

\subsection*{Acknowledgments}
This work was initiated while the second named author was a visiting researcher at 
Saarland University. He would like to express sincere gratitude to the university 
for its hospitality and for providing a perfect research environment.
We thank Kiran Kedlaya for a friendly exchange on pointless K3 surfaces and Asher Auel, Kiran Kedlaya, Philip Engel and Sair Shaikh for providing an advance copy of \cite{AEKS}. 
The first named author gratefully acknowledges support from the Deutsche Forschungsgemeinschaft --- Project-ID 286237555 --- TRR 195.
The second named author is supported by JSPS KAKENHI Grant Number JP24KJ0044. 
The second named author thanks OpenAI for providing complimentary access to 
ChatGPT through the ChatGPT for Academic Researchers program.

\section{Taelman's theorem}\label{sec:Taelman's-theorem}
Let $q$ be a power of a prime $p$, and set $k \coloneq \bF_q$. 
Let $X$ be an ordinary K3 surface over $k$, and fix an embedding $\iota\colon W(k)\hookrightarrow\bC$. As recalled in the Introduction, the canonical lift of $X$ and the embedding $\iota$ give a complex K3 surface $X_{\textrm{can}}^\iota$, and the geometric Frobenius of $X$ induces an endomorphism $F$ of the lattice 
$\Lambda \coloneq H^2(X_{\textrm{can}}^\iota, \bZ)$. 
Here, a \emph{lattice} is a finitely generated free $\bZ$-module equipped with a nondegenerate symmetric bilinear form, and $\Lambda$ is equipped with the intersection form $b:\Lambda\times\Lambda\to \bZ$. 
Note that 
\begin{equation}\label{M1}
\text{$\Lambda$ is a \emph{K3 lattice}, i.e., 
an even unimodular lattice of signature $(3, 19)$.}  \tag{M1} 
\end{equation}
The endomorphism $F$ is a similitude with multiplier $q^2$, that is,  
\begin{equation}\label{M2}
\text{$b(Fx, Fy) = q^2 b(x,y)$ for any $x,y\in \Lambda$.}  \tag{M2} 
\end{equation}
Throughout this paper, the term similitude always means a similitude with multiplier $q^2$. 
We refer to a pair consisting of a lattice and its similitude as a \emph{lattice with similitude}. 
The lattice with similitude $(\Lambda, F)$ satisfies
\begin{equation}\label{M3}
\begin{minipage}{13.5cm}
$F:\Lambda\otimes\bC \to \Lambda\otimes\bC$ is semi-simple and all its eigenvalues have 
absolute value $q$,  
\end{minipage}\tag{M3}
\end{equation}
and 
\begin{equation}\label{M4}
\begin{minipage}{13.5cm}
the $\bZ_p$-module $\Lambda\otimes\bZ_p$ decomposes as 
$\Lambda^0\oplus \Lambda^1\oplus \Lambda^2$ with 
\begin{itemize}
\item[(a)] $F\Lambda^s = q^s \Lambda^s$ for $s=0,1,2$,  
\item[(b)] $\Lambda^0$, $\Lambda^1$, and $\Lambda^2$ are free $\bZ_p$-modules of rank
$1$, $\rk \Lambda - 2$, and $1$, respectively. 
\end{itemize}
\end{minipage}\tag{M4}
\end{equation}

For any lattice with similitude $(\Lambda, F)$ satisfying the conditions \eqref{M1}--\eqref{M4}, we define 
\[ \NS(\Lambda,F) 
\coloneqq \{ x\in \Lambda \mid F^dx = q^dx \text{ for some $d\geq 1$} \}. 
\]
One readily checks that $\NS(\Lambda,F)$ is a primitive $F$-invariant submodule of $\Lambda$. 
Let $\cL$ denote the groupoid whose objects are all lattices with similitude $(\Lambda,F)$ 
with the conditions \eqref{M1}--\eqref{M4} and
\begin{equation}\label{M4half}
\text{$\NS(\Lambda, F)$ is of signature $(1, \rk \NS(\Lambda, F)-1)$.}  \tag{{M4$\unitfrac{1}{2}$}} 
\end{equation}
A morphism from $(\Lambda, F)$ to $(\Lambda', F')$ in $\cL$ is an isometry of 
lattices with similitude, that is, an isometry $t:\Lambda\to \Lambda'$ 
such that $t\circ F = F'\circ t$. 
For $(\Lambda,F)\in \cL$, a connected component of 
\[ 
\left\{ x\in \NS(\Lambda,F)_\bR \;\middle|\; 
\begin{array}{l}
\text{$b(x,x) > 0$, and  } \\
\text{$b(x, \delta)\neq 0$ for any $\delta\in \NS(\Lambda,F)$ with 
$b(\delta,\delta)=-2$}
\end{array}
\right\}
\]
is called a \emph{Weyl chamber}. 
Let $\cT$ denote the groupoid whose objects are all triples $(\Lambda,F,\cK)$, where 
$(\Lambda,F)\in \cL$, and $\cK$ is a Weyl chamber, such that 
\begin{equation}\label{M5}
F\cK = \cK. \tag{M5}
\end{equation}
A morphism from $(\Lambda,F,\cK)$ to $(\Lambda',F',\cK')$ in $\cT$ is 
a morphism $t:(\Lambda,F)\to (\Lambda',F')$ in $\cL$ such that $t\cK=\cK'$. 

We now return to the geometric situation.
Let $X$ be an ordinary K3 surface over $k$, and fix an embedding $W(k) \to \bC$. 
The construction above gives an object $(\Lambda, F)$ of $\cL$. 
There is a natural bijection $\NS(X_{\overline{k}}) \to \NS(\Lambda, F)$.  
Let $\cK\subset \NS(\Lambda, F)_\bR\coloneq \NS(\Lambda, F)\otimes\bR$ 
be the image of the ample cone of $X_{\overline{k}}$
under $\NS(X_{\overline{k}})_\bR \to \NS(\Lambda, F)_\bR$. 
Then, the triple $(\Lambda, F, \cK)$ is an object of $\cT$.   
The assignment $X\mapsto (\Lambda, F, \cK)$ gives rise to a functor $\Phi$
from the groupoid of ordinary K3 surfaces over $k$ to $\cT$. 

\begin{definition}[{\cite{Liedtke-Matsumoto2018}}]
Let $\cO_K$ be a complete discrete valuation ring with field of fractions $K$. 
A K3 surface $X$ over $K$ \emph{satisfies \hypertarget{ast}{\condAst}} if there exists a finite extension 
$L/K$ and a regular algebraic space $\mathfrak{X}$ over $\cO_L$ such that 
$\mathfrak{X}_L \cong X_L$, the special fiber is a strict normal crossing divisor in 
$\mathfrak{X}$, and the relative dualizing sheaf of $\mathfrak{X}/ \cO_L$ is trivial. 
\end{definition}

\begin{theorem}[{\cite[Theorem B]{Taelman20}, \cite{Nygaard1983,Yu2012}}]
\label{th:Taelman2020}
The functor $\Phi$ is fully faithful. If every K3 surface over the field of fractions 
of $W(k)$ satisfies \condAst, then $\Phi$ is essentially surjective. 
\end{theorem}

It is not known whether \condAst always holds. The proof of \Cref{th:has-rational-pts} only needs the existence of the functor $\Phi$ and is independent of \condAst.

\begin{remark}
Conditions \eqref{M1}--\eqref{M4} do not imply the hyperbolicity condition \eqref{M4half}.  
In \cite{Taelman20}, this additional condition appears to be assumed implicitly, 
although it is not stated explicitly. 
\end{remark}

Let $X$ be an ordinary K3 surface over $k$, and write $\Phi(X) = (\Lambda, F, \cK)$. 
For every prime $\ell\neq p$, the comparison isomorphism identifies the action of $F$ on $\Lambda\otimes\bQ_\ell$ with the geometric Frobenius on 
$H^2_{\textrm{\'et}}(X_{\overline{k}},\bQ_\ell)$; see \cite[Theorem 3.1]{Taelman20}. 
Hence, the trace formula together with the vanishing of the odd cohomology yields 
\begin{equation}\label{eq:point-count}
\#X(k) = 1 + \Tr(F) + q^2. 
\end{equation} 

For $(\Lambda, F)\in \cL$, we define the \emph{transcendental lattice} of $(\Lambda,F)$ by
\[ T(\Lambda, F) \coloneq \NS(\Lambda, F)^\perp \subset \Lambda. \] 
The following lemma is obvious but useful.  

\begin{lemma}\label{lem:mor-in-cL}
Let $(\Lambda, F), (\Lambda', F')\in \cL$. An isometry $t:\Lambda\to \Lambda'$
is a morphism in $\cL$ if and only if 
\begin{enumerate}
\item $t(\NS(\Lambda, F)) = \NS(\Lambda', F')$,
or equivalently $t(T(\Lambda, F)) = T(\Lambda', F')$; and 
\item $t|_{\NS(\Lambda, F)}\circ F|_{\NS(\Lambda, F)} = F'|_{\NS(\Lambda', F')} \circ t|_{\NS(\Lambda, F)}$ 
and 
$t|_{T(\Lambda, F)}\circ F|_{T(\Lambda, F)} = F'|_{T(\Lambda', F')} \circ t|_{T(\Lambda, F)}$.
\end{enumerate}
\end{lemma}

\section{Frobenius lattices}\label{sec:Frob-lattices}
We fix a prime $p$ and a power $q = p^e$, where $e\geq 1$.  
The purpose of this section is to collect some properties of $(\Lambda,F)\in\cL$. 
Many of the properties that we need follow solely from Conditions \eqref{M2}, \eqref{M3}, and \eqref{M4}.
To isolate these consequences, we introduce the following notion of a Frobenius lattice. 

\begin{definition}
A \emph{Frobenius lattice} is a lattice with similitude $(L, F)$ satisfying 
Conditions \eqref{M2}, \eqref{M3}, and \eqref{M4}. 
Let $(L,F)$ be a Frobenius lattice.  As in \S\ref{sec:Taelman's-theorem}, define 
\[ \NS(L,F)\defeq \{x \in L \mid F^d x = q^d x\text{ for some }d \geq  1\}.\]
We write $\rho(L,F) \coloneqq \rk\NS(L,F)$ for its rank. 
We say that $(L,F)$ is \emph{effective} if $\NS(L,F)$ has hyperbolic signature, 
and $F$ preserves a Weyl chamber of $\NS(L,F)$. 
We say that $(L,F)$ is \emph{of K3 type} if $L$ is a K3 lattice and $\NS(L,F)$ has 
hyperbolic signature. 
This is equivalent to requiring \eqref{M1} and \eqref{M4half}, that is, $(L,F)\in\cL$.
If $(L,F)$ is effective and of K3 type, then $(L,F,\cK)\in \cT$ where $\cK$ is a Weyl chamber 
preserved by $F$. 
\end{definition}
  

For a lattice with similitude $(L,F)$ (with multiplier $q^2$), it is convenient to 
consider the endomorphism 
\[ q^{-1}F:L_\bQ \to L_\bQ \]
because it is an isometry of $L_\bQ = L\otimes \bQ$, although it need not preserve the lattice $L$. 
We use the following terminology to describe properties of characteristic polynomials of isometries.  

\begin{definition}
Let $K$ be a field, and let $P(X)\in K[X]$ be a monic polynomial. 
We say that $P(X)$ is \emph{$+1$-symmetric} (resp. \emph{$-1$-symmetric}) if 
$P(X) = X^{\deg P} P(X^{-1})$ (resp. $P(X) = -X^{\deg P} P(X^{-1})$). 
If $P$ is either $+1$-symmetric or $-1$-symmetric, then we say that it is \emph{$*$-symmetric}.  
\end{definition}

The characteristic polynomial of an isometry of an inner product space, 
i.e., a finite-dimensional vector space equipped with a nondegenerate symmetric bilinear form, 
is $*$-symmetric. 

In the remainder of this section, we fix a Frobenius lattice $(L,F)$, and write
$\NS = \NS(L, F)$ and $\rho = \rho(L,F)$ briefly. Put 
\[ f \coloneq q^{-1}F:L_\bQ\to L_\bQ, \] 
and let $P(X)\in \bQ[X]$ be the characteristic polynomial of $f$. 
Since $f$ is an isometry, $P(X)$ is $*$-symmetric.
Furthermore $f:L_\bC\to L_\bC$ is semi-simple and all its eigenvalues have absolute 
value $1$ by \eqref{M3}. 
Let $C(X)\in \bZ[X]$ be the largest factor of $P(X)$ that is a product of cyclotomic polynomials.  
Then $P(X)$ decomposes as 
\[
P(X) = S(X) C(X)
\]
where $S(X)\in \bQ[X]$ is the remaining factor. 

\begin{proposition}\label{prop:S-is-*-symm}
We have
$S(1)S(-1)\neq 0$, and $S(X)$ is $+1$-symmetric. In particular, 
the degree of $S(X)$ is even.
\end{proposition}
\begin{proof}
By the definition of $C(X)$, 
we have $(X-1)\nmid S(X)$ and $(X+1)\nmid S(X)$. These show that $S(1)S(-1)\neq 0$. 
Moreover, since $P(X)$ and $C(X)$ are $*$-symmetric, so is $S(X)$. Hence,  
$S(X)$ is a $+1$-symmetric polynomial; because otherwise $S(1) = -S(1)$, 
which would give $S(1) = 0$. 
Then it follows by $S(-1)\neq 0$ that $\deg S$ is even. 
\end{proof}

For a $p$-adic number $\lambda$, let $v_p(\lambda)$ denote the $p$-adic valuation of $\lambda$. 
Note that $v_p$ uniquely extends to the algebraic closure $\overline{\bQ_p}$.
Recall that we have the decomposition $L\otimes \bZ_p = L^0\oplus L^1 \oplus L^2$ 
by \eqref{M4}.  
Let $\lambda_0$ and $\lambda_2\in \bQ_p$ be eigenvalues of $f:L_{\bQ_p}\to L_{\bQ_p}$
corresponding to the eigenlines $L^0_{\bQ_p}$ and $L^2_{\bQ_p} \subset L_{\bQ_p}$, 
respectively. Then 
\[ v_p(\lambda_0) = v_p(q^{-1}) = -e, \quad v_p(\lambda_2) = v_p(q) = e. \]

\begin{proposition}\label{prop:NS}
The lattice $\NS \subset L$ is the largest $f$-invariant sublattice  
on which $f$ acts with finite order. Moreover, we have 
\[ \NS = \{ x\in L \mid C(f)x=0 \}. \]
\end{proposition}
\begin{proof}
By definition, we have  
\begin{equation}\label{eq:NS}
\NS = \{ x\in L \mid \text{$f^d x = x$ for some $d\geq 1$} \}. \tag{$*$}
\end{equation}
We first show that $f$ preserves $\NS$.  
It suffices to show that $f(\NS_\ell) \subset \NS_\ell$
for every prime $\ell$, where $\NS_\ell \coloneq \NS\otimes\bZ_\ell$. Note that 
$\NS_{\bQ_\ell} = \{ x\in L_{\bQ_\ell} \mid \text{$f^d x = x$ for some $d\geq 1$} \}$ 
is preserved by $f:L_{\bQ_\ell} \to L_{\bQ_\ell}$ for every $\ell$. 
For $\ell\neq p$, since $q^{-1}$ is a unit of $\bZ_\ell$, 
it follows that $f$ preserves $L_\ell \coloneq L\otimes \bZ_\ell$ and hence 
$\NS_\ell = \NS_{\bQ_\ell} \cap L_\ell$. 
At $p$, we have $v_p(\lambda_0) = -e <0$ and $v_p(\lambda_2) = e > 0$, whereas  
every root of unity has $p$-adic valuation $0$. 
This implies that $\NS_{\bQ_p}$ is contained in $L^1_{\bQ_p}$, 
and hence $\NS_p = \NS_{\bQ_p}\cap L^1$. 
This intersection is preserved by $f$ since $f(L^1) = L^1$. 
Therefore $f$ preserves $\NS$.  

We next show that $\NS$ is the largest sublattice of $L$ 
on which $f$ acts with finite order.
It is clear by \eqref{eq:NS} that any $f$-invariant sublattice 
on which $f$ acts with finite order is contained in $\NS$.  
Furthermore, the order of $f|_{\NS}:\NS \to \NS$ is 
finite since $\NS$ has finite rank. 

Finally, by the maximality of $C(X)$ and the semi-simplicity of $f$, 
the subspace $\{ x\in L_\bQ \mid C(f)x=0 \}$ is
the largest $f$-invariant subspace of $L_\bQ$ on which $f$ acts with finite order. 
Since $\NS_\bQ$ has the same property, we have
$\NS_\bQ = \{ x\in L_\bQ \mid C(f)x=0 \}$. 
By intersecting with $L$, we obtain 
$\NS = \NS_\bQ \cap L = \{ x\in L \mid C(f)x=0 \}$. 
\end{proof}


\begin{corollary}\label{cor:rho}
We have $\rho = \deg C$. In particular, $\rho\equiv \rk L \bmod 2$. 
\end{corollary}
\begin{proof}
We have $\rho = \rk \NS = \deg C$ by \cref{prop:NS}. Thus 
\[ \rho = \rk L - \deg S \equiv \rk L \mod 2 \]
by \cref{prop:S-is-*-symm}. 
\end{proof}

\begin{corollary}\label{cor:C(1)=0}
Suppose that $(L,F)$ is effective. Then $C(1) = 0$. 
\end{corollary}
\begin{proof}
Let $d$ be the order of $f|_{\NS}$, which is finite by \cref{prop:NS}. 
Pick a nonzero vector $v\in \NS$ that lies in a Weyl chamber $\cK$ preserved by $f$. 
Since $\cK$ is a convex cone, the sum $u \coloneq \sum_{i=0}^{d-1}f^i(v)$ lies in the Weyl chamber, 
and in particular, $u\neq 0$. 
Moreover $f(u) = u$, which shows that $1$ is an eigenvalue of $f$. 
\end{proof}

Since every root of unity has $p$-adic valuation $0$, 
the eigenvalues $\lambda_0$ and $\lambda_2$ of $f:L_{\bQ_p}\to L_{\bQ_p}$ 
are roots of $S(X)$. Consequently, $S(X)\neq 1$, $\NS \subsetneq L$, and $\rho < \rk L$. 
For example, if $(L, F)$ is of K3 type, 
then $\rho$ is an even integer at most $20$ by \cref{cor:rho}. 

In the rest of this section, we collect properties of the factor $S(X)$. 

\begin{proposition}\label{prop:integralcoeff}
We have $P(X) \notin \bZ[X]$ but $qP(X) \in \bZ[X]$, or equivalently, 
$S(X) \notin \bZ[X]$ but $qS(X) \in \bZ[X]$. 
\end{proposition}
\begin{proof}
Let $P^1(X)\in \bZ_p[X]$ be the characteristic polynomial of $f|_{L^1}:L^1\to L^1$. 
Then $P(X)$ decomposes as 
\[P(X) = P^1(X)(X-\lambda_0)(X-\lambda_2)\] 
in $\bQ_p[X]$. 
This shows that $P(X)\notin \bZ_p[X]$ but $qP(X)\in \bZ_p[X]$. 
On the other hand, the coefficients of $P(X)$ belong to $\bZ[p^{-1}]\subset \bQ$
since $F = qf$ preserves the $\bZ$-module $L$. 
Using $\bZ[p^{-1}] \cap \bZ_p = \bZ$, we obtain  
$P(X)\notin \bZ[X]$ but $qP(X)\in \bZ[X]$. 
This is equivalent, by Gauss's lemma, to the condition that 
$S(X) \notin \bZ[X]$ but $qS(X) \in \bZ[X]$. 
\end{proof}

\begin{proposition}\label{prop:S-is-irr}
The factor $S(X)$ is irreducible over $\bQ$. 
\end{proposition}
\begin{proof}
Let $S'(X)\in \bQ[X]$ be the irreducible monic factor of $P(X)$ over $\bQ$ such that 
$X-\lambda_0$ divides $S'(X)$ over $\bQ_p$. Then, $S'(X)$ divides $S(X)$ over $\bQ$ 
since $\lambda_0$ is not a root of unity, and we can write 
$S(X) = S'(X)Q(X)$ for some monic polynomial $Q(X)\in \bQ[X]$. 
The same argument as in \cref{prop:integralcoeff} shows that $qS'(X)\in \bZ[X]$ 
and that this polynomial is primitive. 
Since $Q(X)$ is the quotient of $qS(X)\in\bZ[X]$ by $qS'(X)\in \bZ[X]$, 
Gauss's lemma implies that $Q(X)\in \bZ[X]$.  
On the other hand, since all roots of $Q(X)$ have absolute value $1$, 
$Q(X)$ is a product of cyclotomic polynomials by Kronecker's theorem. 
Therefore $Q(X)$ must be the constant $1$ by the maximality of $C(X)$. 
This shows that $S(X) = S'(X)$, which is irreducible over $\bQ$.  
\end{proof}

\begin{proposition}\label{prop:S3}
The eigenvalue $\lambda_2\in \bQ_p$ is the unique root of $S(X)$ with 
$p$-adic valuation $e$, and all roots in $\overline{\bQ_p}$ 
except for $\lambda_0$ and $\lambda_2$ have $p$-adic valuation $0$. 
Furthermore, we have $\lambda_0 = \lambda_2^{-1}$. 
\end{proposition}
\begin{proof}
In $\bQ_p[X]$, we can write $S(X) = S^1(X)(X-\lambda_0)(X-\lambda_2)$ for some $S^1(X)\in \bQ_p[X]$. 
Since the roots of $S^1(X)$ are among the eigenvalues of the isometry $f|_{L^1}$ 
of the $\bZ_p$-lattice $L^1$, all roots of $S^1(X)$ have $p$-adic valuation $0$. 
Thus $\lambda_0$ and $\lambda_2$ are the unique roots of $S$ with valuations $-e$ and $e$, respectively.
Since $S$ is $+1$-symmetric, every root occurs together with its inverse. 
Therefore $\lambda_0 = \lambda_2^{-1}$, and this completes the proof. 
\end{proof}

To summarize the properties of $S(X)$, we introduce $q$-fractional Salem polynomials, 
which can be viewed as an analogue of Salem polynomials. 

\begin{definition}
A \textit{$q$-fractional Salem polynomial} is a polynomial $S(X) \in \bQ[X]$ satisfying the following conditions:
\begin{align}
\label{S1}
&\text{$S(X)$ is an irreducible monic polynomial, and $qS(X)\in \bZ[X]$.} \tag{S1}\\
\label{S2}
&\text{All roots in $\bC$ of $S(X)$ have absolute value $1$. } \tag{S2}\\
\label{S3}
&\begin{minipage}{11.5cm}
There exists a unique root $\lambda\in \bQ_p$ of $S(X)$ such that 
$v_p(\lambda) = e$, and all roots in $\overline{\bQ_p}$ 
except for $\lambda$ and  $\lambda^{-1}$ have $p$-adic valuation $0$. 
\end{minipage}
\tag{S3}
\end{align}
\end{definition}

\begin{corollary}
Let $(\Lambda, F)$ be a Frobenius lattice, and let $S(X)$ be the factor of the characteristic polynomial 
of $q^{-1}F$ defined as above. Then $S(X)$ is a $q$-fractional Salem polynomial. 
\end{corollary}
\begin{proof}
Condition \eqref{S1} follows from \cref{prop:integralcoeff,prop:S-is-irr};   
Condition \eqref{S2} follows from the assumption \eqref{M3}; and 
Condition \eqref{S3} follows from \cref{prop:S3}. 
\end{proof}

\section{\texorpdfstring{$q$}{q}-fractional Salem polynomials}\label{sec:qSalem}
As in the previous section, we fix a prime $p$ and a power $q = p^e$, where $e\geq 1$.  
In this section, we give an algorithm for enumerating all $q$-fractional Salem polynomials 
of a prescribed degree. 

Note that any irreducible monic polynomial in $\bQ[X]$ all of whose roots 
have absolute value $1$ is $+1$-symmetric and has even degree, except for
$X-1$ and $X+1$. Since $X-1$ and $X+1$ do not satisfy \eqref{S3}, 
every $q$-fractional Salem polynomial is $+1$-symmetric and has even degree. 
In particular, if $\alpha$ is a root of a $q$-fractional Salem polynomial, then 
$\alpha^{-1}$ is also a root of the polynomial.  
The following proposition gives a coefficient criterion for Condition \eqref{S3}.

\begin{proposition}
Let $S(X) = X^n + a_1 X^{n-1} + \cdots \in \bQ[X]$ be a polynomial satisfying 
\eqref{S1} and \eqref{S2}. 
The following three conditions are equivalent: 
\begin{itemize}
\item[\textup{(a)}] $S(X)$ is a $q$-fractional Salem polynomial, i.e., $S(X)$ satisfies \eqref{S3}. 
\item[\textup{(b)}] There exists a unique root $\lambda\in \overline{\bQ_p}$ of $S(X)$ such that 
$v_p(\lambda) = e$, and all roots in $\overline{\bQ_p}$ 
except for $\lambda$ and $\lambda^{-1}$ have $p$-adic valuation $0$. 
\item[\textup{(c)}] $p^{e-1}a_1 \notin \bZ$. 
\end{itemize}
\end{proposition}
\begin{proof}
(a) $\Rightarrow$ (b). Clear.   

(b) $\Rightarrow$ (c). 
Let $\alpha_1, \ldots, \alpha_n$ be the roots of $S(X)$. 
Since $\lambda^{-1}$ is the unique root of valuation $-e$, the ultrametric inequality gives
\[v_p(a_1) = v_p\left(-\sum_{i=1}^n \alpha_i \right) = v_p(\lambda^{-1}) = -e. \]
Hence, $v_p(p^{e-1}a_1) = -1$, and $p^{e-1}a_1\notin\bZ$. 

(c) $\Rightarrow$ (a). 
Consider the Newton polygon of $S(X)$ with respect to $p$. 
Note that $S(X)$ is $+1$-symmetric. Conditions \eqref{S1} and (c) imply that 
the slopes of the Newton polygon are $-e, 0, e$, and the corresponding lengths are 
$1, n-2, 1$. 
Then, (a) follows from basic results on Newton polygons; see \cite[\S 7.4]{Gouvea_padic}. 
%
\end{proof}

Another important observation is as follows: 

\begin{proposition}\label{prop:finiteness-of-q-Salem}
For each $d \in \bZ_{\geq 1}$, there are only finitely many $q$-fractional Salem polynomials of degree $2d$.
\end{proposition}
\begin{proof}
Let $S(X)$ be a $q$-fractional Salem polynomial of degree $2d$, say 
\[ S(X) = X^{2d} + a_1 X^{2d-1} +\cdots + a_dX^d + \cdots + a_1 X + 1. \]
From Vieta's formulas with the condition \eqref{S2}, it follows that 
$|a_i| \leq \binom{2d}{i}$ for $i = 1, \ldots, d$. Hence, noting that $qS(X)\in \bZ[X]$, 
there are only finitely many $q$-fractional Salem polynomials of fixed degree $2d$. 
\end{proof}

We now describe an algorithm for enumerating $q$-fractional Salem polynomials. 
Our approach is inspired by \cite{mossinghoff:minimal-mahler}. 
A more general search algorithm for root-unitary polynomials is formulated by 
Kedlaya and Sutherland \cite{Kedlaya-Sutherland2016}, 
based on FLINT, which is available e.g. in SageMath. 

Vieta's bounds used in \cref{prop:finiteness-of-q-Salem} for the coefficients of $S(X)$ are too loose for the enumeration of $q$-fractional Salem polynomials in practice. We turn to Newton's identities to reduce the computational cost. 


\begin{definition}
Let 
\[P(X) = X^n + a_1X^{n-1} + a_2 X^{n-2} + \cdots  + a_{n-1} X + a_n \in \bC[X]\]
be a monic polynomial with roots $\alpha_1, \ldots, \alpha_n\in\bC$. 
For $i\geq 1$, put 
\[p_i = \sum_{j=1}^n \alpha_j^i.\] 
The following identities are called \textit{Newton's identities}: 
\[ p_i + p_{i-1}a_{1} + p_{i-2}a_{2} + \cdots + p_{1}a_{i-1} + i a_i = 0 
\quad(1 \leq i \leq n).  \]
Hence, $p_i$ depends only on $a_1,\dots, a_{i}$ and can be computed efficiently.   
\end{definition}
Let $\cS_{q,2d}$ be the set of all $q$-fractional Salem polynomials of degree $2d$, and   
let 
\[S(X) = X^{2d} + a_1X^{2d-1} + \cdots + a_dX^d + \cdots + a_1 X + 1\in \cS_{q,2d}\] 
with roots $\alpha_1, \ldots, \alpha_{2d}\in \bC$. 
In order to enumerate the set $\cS_{q,2d}$, we make the following observations:
\begin{enumerate}
\item We have $p^{e-1}a_1\notin\bZ$. 
Since $S(-X) = X^{2d} - a_1 X^{2d-1} + \cdots - a_1 X + 1$ also belongs to $\cS_{q,2d}$, 
it is enough to list all $S(X)$ with $a_1 > 0$. 
\item Put $p_i = \sum_{j=1}^{2d} \alpha_j^i$ for $i\geq 1$. 
We have $|p_i|<2d$ for $i \geq 1$ by \eqref{S2}. 
\item Let $T(Y)$ be the trace polynomial of $S(X)$, characterized by $S(X)=X^dT(X+X^{-1})$, 
and write $T(Y)=Y^d + b_1 Y^{d-1} + \cdots + b_d$.
All roots of $T(Y)$ are simple, real, and lie in $(-2,2)$. 
Hence, all roots of the $k$-th derivative $T^{(k)}(Y)$ are also simple, real, and lie in $(-2,2)$. 
Note that $T^{(k)}(Y)$ is determined by $a_1, \ldots, a_{d-k}$ since so are 
$b_1, \ldots, b_{d-k}$. 
\end{enumerate}

The enumeration algorithm proceeds by a depth-first tree traversal. The tree structure 
is implicit and only the path from the root to the current vertex is stored at any moment.
The tree has $d+1$ levels. The root of the tree is $a_0=1$ at level $0$.
For prescribed bounds $0\leq l_1<u_1\leq2d$, the vertices at level $1$ are
\[
\left\{ a_1 \in q^{-1}\mathbb{Z} \;\middle|\; 
a_1 \notin p^{-(e-1)}\mathbb{Z} \text{ and } l_1 < a_1 < u_1 \right\}. 
\] 
These bounds allow one to restrict the enumeration to a prescribed range of $a_1$, 
or equivalently of the trace of $S$.
At level $i$, for $2\leq i \leq d$, we have the candidates for $a_i\in q^{-1}\bZ$ which satisfy
\begin{enumerate} 
\item[(ii)] $|p_{i}| < 2d$;
\item[(iii)] all roots of $T^{(d-i)}(Y)$ are simple, real, and lie in $(-2,2)$.
\end{enumerate}

We note that these conditions depend only on $a_0,\dots, a_{i}$.
At level $d$, we also impose 
\begin{enumerate} 
\item[(iv)] the resulting polynomial $S(X)$ is irreducible. 
\end{enumerate}
Testing whether all roots of $T^{(d-i)}(Y)$ are simple, real, and in $(-2,2)$ is possible without actually computing them, namely by Sturm's theorem (see e.g. \cite[Exercise 4.32]{modern-computer-algebra}).

One of the aims of this paper is to search for pointless K3 surfaces. 
Let $(\Lambda, F)$ be an effective Frobenius lattice of K3 type, 
and write $P(X) = S(X)C(X)$ for the decomposition of the characteristic polynomial 
of $q^{-1}F$ from \S \ref{sec:Frob-lattices}. 
Put $d =\deg(S)/2$. By \cref{cor:C(1)=0}, the polynomial $C(X)$ is divisible by $X-1$. 
Since $\deg C=22-2d$ and all roots of $C(X)$ have absolute value $1$, we have
\[
\Tr(C)\geq 1-(21-2d)=2d-20,
\]
where $\Tr(C)$ denotes the sum of the roots of $C(X)$. 
If the K3 surface corresponding to $(\Lambda, F)$ is pointless, then 
\eqref{eq:point-count} and $\Tr(F) = q(\Tr(S) + \Tr(C))$
give the following inequality:  
\begin{equation}\label{eq:trace-bound}
1 + q(2d - 20) + q\Tr(S) + q^2 \leq 0.
\end{equation}

\begin{proposition}
Let $d$ be an integer with $1\leq d \leq 10$, and suppose that $q\leq 19$.  
The number of $q$-fractional Salem polynomials of degree $2d$ 
satisfying the trace bound \eqref{eq:trace-bound}
is given in \cref{table:q-salem-small-trace}. 
\end{proposition}

\begin{table}[h!]
\caption{Numbers of $q$-fractional Salem polynomials satisfying \eqref{eq:trace-bound}\\
\centering{}}\label{table:q-salem-small-trace}
\begin{tabular}{c|cccccccccc}
$q \setminus 2d$ &2 & 4 & 6 &8 &10 &12 & 14 & 16& 18 & 20\\
\hline
2 & 4 & 8 & 30 & 98 & 274 & 792 & 1700 & 4427 & 4601 & 356\\
3 & 8 & 36 & 208 & 1076 & 4828 & 20602 & 76733 & 235299 & 179551 & 6568\\
4 & 8 & 56 & 446 & 3118 & 19156 & 109573 & 541858 & 1411150 & 552590 & 5266\\
5 & 16 & 152 & 1572 & 13652 & 106358 & 760327 & 3986405 & 6371011 & 908997 & 784\\
7 & 24 & 356 & 5146 & 63714 & 699300 & 5625073 & 14006291 & 3865658 & 20806\\
8 & 16 & 280 & 4622 & 65757 & 788865 & 4831810 & 5731839 & 457293 & 61\\
9 & 24 & 484 & 9018 & 143731 & 1623536 & 5814915 & 2593258 & 36238\\
11 & 40 & 1020 & 23185 & 361080 & 1747798 & 1074120 & 23237\\
13 & 48 & 1465 & 31549 & 202301 & 165725 & 5340\\
16 & 32 & 620 & 2293 & 713\\
17 & 48 & 541 & 779 & 28\\
19 & 18 & 20
\end{tabular}
\end{table}

\section{Strategy for enumeration}\label{sec:strategy}
As in the previous section, we fix a prime $p$ and a power $q$ of $p$. 
In this section, we describe the strategy for enumerating the objects of $\cT$, 
that is, effective Frobenius lattices of K3 type. 

We first recall the following constraint on the characteristic polynomial of 
the Frobenius. 

\begin{proposition}[cf. {\cite[Theorem 2]{Kedlaya-Sutherland2016}}]\label{prop:sq}
Let $X$ be an ordinary K3 surface over $k = \bF_q$ and  
$F: \Lambda \to \Lambda$ be the endomorphism of a K3 lattice $\Lambda$ 
described in \S\ref{sec:Taelman's-theorem}. 
Let $P$ be the characteristic polynomial of $q^{-1}F$. 
Then $qP(-1)$ is a square number (possibly $0$). 
\end{proposition}

This proposition follows from the Artin--Tate formula, as explained in 
\cite{Kedlaya-Sutherland2016}.  
It can also be proved by a purely lattice-theoretic argument. 
In fact, it holds for any Frobenius even unimodular lattice.  

We next describe the compatibility between the transcendental part and the 
N\'eron--Severi part. 
Let $(\Lambda, F)\in \cL$ be a Frobenius lattice of K3 type. 
The transcendental lattice $T \coloneq T(\Lambda, F)$ equipped 
with the restriction $F|_T:T\to T$ is a lattice with similitude. 
Moreover, it has the following property. 

\begin{lemma}\label{lem:unimodularity}
The transcendental part $(T,F|_T)$ of $(\Lambda, F)$ is a Frobenius lattice 
such that $T^0 = \Lambda^0$ and $T^2 = \Lambda^2$. Furthermore, 
$T^0\oplus T^2$ is unimodular. 
\end{lemma}
\begin{proof}
Conditions \eqref{M2} and \eqref{M3} for $(T,F|_T)$ follow from the corresponding conditions for $(\Lambda,F)$. 
By considering the decomposition into eigenspaces, one sees that the orthogonal direct summand $\Lambda^0\oplus \Lambda^2$ of $\Lambda_p \coloneq \Lambda\otimes \bZ_p$
is nondegenerate. Since the ambient lattice $\Lambda_p$ is unimodular, it follows that $\Lambda^0\oplus \Lambda^2$ is unimodular. 
Since every root of unity has $p$-adic valuation $0$, we have 
$\NS(\Lambda, F)\otimes \bZ_p  \subset \Lambda^1$. 
Thus $T_p\coloneq T\otimes \bZ_p \supset \Lambda^0\oplus \Lambda^2$. 
Since $\Lambda^0 \oplus \Lambda^2$ is unimodular, this sublattice must be an 
orthogonal direct summand of $T_p$. Namely, we have 
$T_p = \Lambda^0 \oplus T^1 \oplus \Lambda^2$, where 
$T^1$ is the orthogonal complement of $(\Lambda^0\oplus \Lambda^2)$ in $T_p$. 
Putting $T^0 = \Lambda^0$ and $T^2 = \Lambda^2$,
this decomposition gives \eqref{M4} for $(T, F|_T)$, and hence 
$(T, F|_T)$ is a Frobenius lattice. 
Furthermore, $T^0\oplus T^2 = \Lambda^0 \oplus \Lambda^2$ is unimodular.  
\end{proof}

Although the isometry $t\coloneq q^{-1}F|_T:T_\bQ\to T_\bQ$ does not preserve 
the lattice $T$, the unimodularity of $T^0\oplus T^2$ allows us to define 
an action of $t$ on the discriminant group $D_T\coloneq T^\vee/T$ as follows. 
Recall that $D_T$ decomposes as 
\[ D_T = \bigoplus_{\ell\text{: prime}} (D_T)_\ell, \] 
where $(D_T)_\ell$ is the $\ell$-primary component of $D_T$, 
which is canonically identified with 
the discriminant group $D_{T_\ell}$ of $T_\ell \coloneq T\otimes \bZ_\ell$. 
For every prime $\ell \neq p$, the isometry $t:T_{\bQ_\ell}\to T_{\bQ_\ell}$ preserves 
the lattice $T_\ell$ and thus induces an isometry of $D_{T_\ell}=(D_T)_\ell$. 
On the other hand, the unimodularity of $T^0\oplus T^2$ gives a canonical isomorphism 
\[ (D_T)_p = D_{T_p} \cong D_{T^1}. \] 
Since $t:T_{\bQ_p}\to T_{\bQ_p}$ preserves $T^1$,  
its restriction to $T^1$ induces an isometry of $D_{T^1}$, and hence of $(D_T)_p$. 
Combining these local actions, we obtain an isometry of $D_T$, 
which is denoted by $D_t$. 

Let $\NS \coloneq \NS(\Lambda, F)$ and $n\coloneq q^{-1} F|_{\NS}:\NS_\bQ\to \NS_\bQ$. 
By \cref{prop:NS}, $n$ preserves the lattice $\NS$ and thus is an isometry of $\NS$. 
We have a natural \emph{anti-isometry}
\[ \gamma: D_{\NS} \to D_T,  \]
that is, a group isomorphism such that $\bar{x}^2 = -\gamma(\bar{x})^2$
for any $\bar{x}\in D_{\NS}$, 
where, for an even lattice $L$ and an element $z\in L^\vee$, 
$\bar{z}\in D_L$ denotes the class of $z$, and 
$\bar{z}^2$ denotes the value in $\bQ/2\bZ$ of the induced quadratic form. 
As for an orthogonal decomposition of an even unimodular lattice with isometry, 
the induced actions $D_n$ and $D_t$ satisfy the compatibility relation 

\begin{equation}\label{eq:compatibility}
D_t \circ \gamma = \gamma\circ D_n, 
\end{equation}
see e.g. \cite[\S 1.3]{Kondo_K3book}. 

Conversely, let $(L, \hat{t})$ be a Frobenius even lattice  
such that $L^0\oplus L^2$ is unimodular and $\NS(L, \hat{t})=0$, 
and let $(N, n)$ be the pair of an even lattice $N$ and an isometry $n$ of $N$ 
with finite order.  
Note that $(N\oplus L)^\vee = N^\vee \oplus L^\vee$, and thus 
$D_{N\oplus L} = D_N\oplus D_L$. 
Suppose that there exists an anti-isometry $\gamma:D_N\to D_L$ such that
$D_t\circ \gamma = \gamma \circ D_n$, where $t \coloneq q^{-1}\hat{t}$.  
Then there exists a unique even unimodular lattice $\Lambda$ 
with $N \oplus L \subset \Lambda \subset (N \oplus L)^\vee$ 
such that $\Lambda/(N \oplus L)\subset D_{N\oplus L} = D_N\oplus D_L$ 
is the graph of $\gamma$. 
Moreover, 
$F \coloneq q(n\oplus t): N_\bQ \oplus L_\bQ \to N_\bQ \oplus L_\bQ$
preserves the lattice $\Lambda$, and 
the pair $(\Lambda, F)$ is a Frobenius lattice such that 
$(\NS(\Lambda, F), F|_{\NS(\Lambda, F)}) = (N, qn)$ and 
$(T(\Lambda, F), F|_{T(\Lambda, F)}) = (L, \hat{t})$. 
Namely, the anti-isometry $\gamma$ glues $(N, qn)$ and $(L, \hat{t})$ 
together into a Frobenius even unimodular lattice. 

In the algorithm for enumerating Frobenius lattices of K3 type, 
we first list candidates for the transcendental and N\'eron--Severi parts
and then glue them together. 
Thus, we need to determine when two anti-isometries between discriminant groups
give rise to isometric Frobenius lattices.  
Fix an anti-isometry $\gamma:D_N\to D_L$ such that
$D_t\circ \gamma = \gamma \circ D_n$. Then 
the set of equivariant anti-isometries $D_N\to D_L$ can be written as 
\[ \gamma O(D_N, D_n) = \{ \gamma \circ \tau \mid \tau \in O(D_N, D_n) \},   \]
where 
\[ O(D_N, D_n) \coloneq \{ \tau \in O(D_N) \mid \tau \circ D_n = D_n \circ \tau \}. \]
Using \cref{lem:mor-in-cL}, we obtain the following proposition. 

\begin{proposition}
Let $d$ be an integer with $1 \leq d \leq 10$. 
Let $(L, \hat{t})$ be an even Frobenius lattice of signature $(2, 2d-2)$ 
such that $L^0\oplus L^2$ is unimodular and $\NS(L, \hat{t}) = 0$, 
and let $(N, n)$ be a pair consisting of an even lattice $N$ of signature $(1, 21-2d)$
and an isometry $n$ of $N$ with finite order. 
Let $\gamma \colon D_N \to D_L$ be an anti-isometry
with $D_t\circ \gamma=\gamma\circ D_n$. 
Then the set of isomorphism classes of Frobenius lattices of K3 type $(\Lambda, F)$ with 
$(T, F|_T) \cong (L, \hat t)$ and $(\NS, q^{-1}F|_{\NS}) \cong (N,n)$ is in bijection with the double coset 
\begin{equation}\label{eq:double-coset}
O(L, t) \setminus \gamma O(D_N, D_n) / O(N,n).
\end{equation}
\end{proposition}

We now describe our algorithm for enumerating effective Frobenius lattices of K3 type. 

\begin{steps}
\item For each $1\leq d \leq 10$, 
enumerate all $q$-fractional Salem polynomials $S$ of degree $2d$ and 
all products of cyclotomic polynomials $C$ of degree $22-2d$. 
Enumeration of $q$-fractional Salem polynomials of fixed degree is described in 
\S\ref{sec:qSalem}, and it is easy to compute all 
cyclotomic polynomials $C$ of fixed (small) degree. 

\item Fix a $q$-fractional Salem polynomial $S$ of degree $2d$
and a product of cyclotomic polynomials $C$ of degree $22-2d$ such that 
$C(1)=0$ and $qS(-1)C(-1)$ is a square. 
The product $SC$ will be the characteristic polynomial of $q^{-1}F$
for the resulting Frobenius lattice $(\Lambda, F)$ of K3 type. 
The condition $C(1)=0$ is imposed by \cref{cor:C(1)=0}, and 
the condition that $qS(-1)C(-1)$ is a square is imposed by \cref{prop:sq}. 

\item Compile a list of all candidates for transcendental parts $(L, \hat{t})$ 
such that $q^{-1}\hat{t}$ has characteristic polynomial $S$. 
The information in $C$ restricts the possible transcendental parts to 
be a finite set. 
We will explain the details in \S\ref{sec:transcendental-part}.

\item Fix $(L, \hat{t})$ from Step 3. Compile a list of all candidates for  
N\'eron--Severi parts $(N, n)$ compatible with $(L, \hat{t})$. 
Effectivity can be checked at this step. 
We will explain the details in \S\ref{sec:Neron--Severi-part}.  

\item Fix $(N, n)$ from Step 4. Compute the double coset \eqref{eq:double-coset}. 
Using each anti-isometry that represents an element of the double coset, 
glue $(N, n)$ and $(L, \hat{t})$ to obtain a Frobenius lattice of K3 type. 
\end{steps}

\section{The transcendental part}\label{sec:transcendental-part}
As in the previous section, we fix a prime $p$ and a power $q$ of $p$. 
In this section, we explain how to enumerate all candidates for transcendental parts 
of a Frobenius lattice of K3 type whose characteristic polynomial is prescribed, 
beginning with a more general setting. 

\subsection{Lattices with similitude 
whose characteristic polynomials are powers of irreducible polynomials}
In this subsection, we study lattices with similitude whose characteristic polynomial 
is a power of a single irreducible polynomial, and describe a way to enumerate them. 
We fix an irreducible $+1$-symmetric polynomial $S(X)\in\bQ[X]$ of even degree 
such that $qS(X)\in \bZ[X]$ and a positive integer $m$.
We define a field $E$ and its element $\alpha$ by 
\[E \coloneq \bQ[X]/(S(X)) \quad \mbox{ and } \quad \alpha \coloneq X + (S) \in E. \] 
We often use the symbol $\hat{t}$ for a similitude and write $t = q^{-1}\hat{t}$, 
the associated isometry (over $\bQ$).  

\begin{definition}
A \emph{fractional lattice} is a finitely generated free $\bZ$-module equipped with 
a non-degenerate symmetric bilinear form $b:L\times L\to \bQ$. 
We call a usual lattice an \emph{integral} lattice 
when we want to emphasize that its bilinear form takes values in $\bZ$. 
We define $\cF_{S^m}$ to be the groupoid whose objects are fractional 
lattices with similitude $(L, \hat{t})$
such that $\hat{t}$ is semisimple and the characteristic polynomial of 
$t\coloneq q^{-1}\hat{t}: L_\bQ\to L_\bQ$ is equal to $S(X)^m$. 
Although this is an abuse of notation, we often write $(L, t)\in \cF_{S^m}$ instead of 
$(L, \hat{t})\in \cF_{S^m}$. We define 
\[ O(L,t) \coloneq \Hom_{\cF_{S^m}}((L,t), (L,t)) 
= \{ \theta\in O(L) \mid \theta\circ t = t \circ \theta \}. 
\]
\end{definition}

Let $(L, t)\in \cF_{S^m}$, 
and put $V \coloneq L\otimes \bQ$. 
Then $V$ is an $E$-vector space of dimension $m$ by the following action: 
\[ g(X).v = g(t).v \quad (g(X)\in \bQ[X], v\in V).  \]
Let $\hat{\alpha} \coloneq q\alpha \in \cO_E$. By the same action,
the lattice $L$ is a $\o$-module. 
When $\o$ agrees with the maximal order $\O_E$ of $E$, which is a Dedekind domain,
such lattices are easier to handle. For example,  
we can use the theory of hermitian $\O_E$-lattices, as we explain below. 
In general, it is convenient to consider the largest $\cO_E$-module that is contained in $L$, 
denoted $M_L$. 
Note that $(M_L, t)$ (precisely, $(M_L, \hat{t}|_{M_L})$) also belongs to $\cF_{S^m}$.

\begin{proposition}
Let $(L,t)$ and $(L',t')\in \cF_{S^m}$. If $(L, t)\cong (L', t')$ 
then $(M_L, t) \cong (M_{L'} , t')$. 
In other words, $(M_L, t)$ is an invariant for $(L, t)$. 
\end{proposition}
\begin{proof}
Straightforward. 
\end{proof}

To enumerate objects in $\cF_{S^m}$, 
we first enumerate possible invariants $(M, s)$ and then list 
all lattices with similitude $(L, t)$ such that $(M_L, t) \cong (M, s)$. 
In the situation relevant for us, we have a constraint on the discriminant group, 
and it ensures that only finitely many objects need to be considered.  

\begin{definition}
A \emph{discriminant annihilator} for an integral lattice with similitude 
$(L, t)\in \cF_{S^m}$ is an ideal 
$I$ of $\cO_E$ or of $\o$ such that $I L^\vee \subset L$. 
\end{definition}

We will give an example of a discriminant annihilator for the transcendental part 
of a Frobenius lattice of K3 type in the next subsection (\cref{ex:disc-ann}). 
In what follows, the conductor of $\o$, which quantifies the difference between 
$\o$ and $\cO_E$, is useful. 

\begin{definition}
The \emph{conductor} of $\o$, denoted $\mathfrak{f}$, 
is the annihilator of the $\o$-module $\O_E/\o$:
\[\mathfrak{f} \defeq \mathrm{Ann}_\o(\O_E/\o) 
= \{a \in E \mid a \O_E \subset \o \}.\]
It is the largest $\O_E$-module contained in $\o$. 
\end{definition}

For any $(L, t)\in \cF_{S^m}$, we have $\frakf L \subset \o L \subset L$, and hence 
\begin{equation}\label{eq:fL<M_L}
\frakf L \subset M_L
\end{equation}
by the maximality of $M_L$.  
A discriminant annihilator for $(L, t)$ gives one for $(M_L, t)$. 

\begin{proposition}\label{prop:discriminant-annihilator}
Let $(L,t)\in \cF_{S^m}$ be an integral lattice with similitude, 
and let $I$ be a discriminant annihilator for $(L,t)$. 
Then $\frakf I$ is a discriminant annihilator for $(M_L,t)$. 
\end{proposition}
\begin{proof}
By \eqref{eq:fL<M_L}, we have $\cO_E L \subset \frakf^{-1}M_L$. 
By $IL^\vee \subset L$, we also have $I\cO_E L^\vee \subset \cO_E L$. 
On the other hand, we have $\cO_E L^\vee = M_L^\vee$. 
Indeed, $\cO_E L^\vee$ is the smallest $\cO_E$-module containing $L^\vee$, 
which coincides with $M_L^\vee$, since $M_L$ is the largest $\cO_E$-module contained 
in $L$. Thus, we obtain 
\[ \frakf I M_L^\vee = \frakf (I\cO_EL^\vee) \subset \frakf (\cO_E L) 
\subset\frakf (\frakf^{-1}M_L) = M_L. \]
\end{proof}

\subsubsection{Case where the underlying lattice is an $\cO_E$-module}
Here, we describe how to enumerate objects in $\cF_{S^m}$ 
whose underlying lattice is an $\cO_E$-module. 
To this end, we introduce the notion of a hermitian lattice. 

\begin{definition}
Let $E$ be an algebraic number field with an involution $\sigma:E\to E$ and $K$ be the fixed field.
A \emph{fractional hermitian $\cO_E$-lattice}, or just 
\emph{fractional hermitian lattice}, is a finitely generated projective 
$\cO_E$-module $M$ equipped with a non-degenerate hermitian form 
$h:M\times M \to E$. Here, we use the convention that a hermitian form 
is linear in its first variable and satisfies $h(x,y)=\sigma h(y,x)$ for all $x,y\in M$. 
The \emph{scale ideal} $\scale(M) = \scale(M, h)$ of a fractional hermitian lattice $(M, h)$ 
is the fractional ideal of $\cO_E$ generated by the set $\{ h(x,y)\mid x, y\in M \}$. 
A morphism of fractional hermitian lattices is an $\cO_E$-homomorphism 
compatible with the hermitian forms.
\end{definition}

We remark that the field $E = \bQ[X]/(S(X))$ admits an involution $\sigma$
determined by $\alpha\mapsto \alpha^{-1}$. 
For a hermitian form $h$ on the $E$-vector space $V = L\otimes \bQ$, 
the restriction of scalars
\[ \Tr_{E/\bQ}\circ h : V \times V \to \bQ \]
defines a $\bQ$-bilinear form on $V$. 
Conversely, if $b:V\times V\to \bQ$ is a non-degenerate symmetric bilinear form for which $t$ is an isometry, 
then there exists a unique hermitian form $h$ on $V$ satisfying 
$b = \Tr_{E/\bQ}\circ h$; see \cite[Lemma 1.1]{Milnor1969}.  
Therefore, we obtain:

\begin{proposition}
Restriction of scalars $(M, h)\mapsto (M, \Tr_{E/\bQ}\circ h)$, 
with the similitude given by multiplication by $q\alpha$, defines 
an equivalence between the category of fractional hermitian $\cO_E$-lattices 
of rank $m$ and the full subcategory of $\cF_{S^m}$ consisting of objects whose 
underlying lattice is an $\cO_E$-module.  
\end{proposition}

Thanks to Kirschmer \cite{Kirschmer2019, KirschmerHabil}, algorithms for working with hermitian lattices are available. 
These algorithms use the descriptions of their genera due to Jacobowitz \cite{jacobowitz} 
and are implemented, for example, in OSCAR \cite{OSCAR}. 
For example, we can compute all genera and hence all isomorphism classes 
of hermitian lattices $(M, h)$ with the bounds $\scale(M) \subset \mathfrak{a}$
and $\scale(M^\#) \subset \mathfrak{b}$ for prescribed fractional ideals 
$\mathfrak{a}$ and $\mathfrak{b}$ of $\cO_E$, where 
\[M^\# \coloneq \{x \in M \otimes_{\cO_E} E \mid h(x,y) \in \O_E 
\text{ for all $y \in M$} \} \]
is the hermitian dual. 
Such bounds can be obtained from a discriminant annihilator. 

\begin{proposition}
Let $(M,t)\in \cF_{S^m}$ be an integral lattice with similitude, and 
let $h$ denote the hermitian form associated with the bilinear form on $M$. 
Let $\frakD \coloneq \frakD_{E/\bQ}$ be the different ideal of $E/\bQ$. 
Then $\scale(M, h) \subset \frakD^{-1}$. 
If $J$ is a discriminant annihilator for $(M,t)$, then 
 $\scale(M^\#, h) \subset \frakD J^{-1} \sigma(J^{-1})$.   
\end{proposition}
\begin{proof}
Note that $\frakD^{-1} = \{ a\in E \mid \Tr_{E/\bQ}(ab)\in\bZ \text{ for all $b\in \cO_E$} \}$.
Since $M = (M, b)$ is an integral lattice, the relation $b = \Tr_{E/\bQ}\circ h$ 
shows that $\scale(M, h) \subset \frakD^{-1}$. 

Suppose that $J$ is a discriminant annihilator for $(M,t)$, i.e., $JM^\vee \subset M$.  
One easily sees that $\frakD^{-1}M^\# \subset M^\vee$. Hence 
$M^\# \subset \frakD M^\vee \subset \frakD J^{-1} M$.
Taking the scale ideals, we obtain 
\[ \scale(M^\#) \subset \scale(\frakD J^{-1} M) 
= \frakD J^{-1} \sigma(\frakD J^{-1})\scale(M)
\subset \frakD J^{-1} \sigma(J^{-1}),  
\]
where, in the last inclusion, we use $\sigma(\frakD) = \frakD$ and 
$\scale(M, h) \subset \frakD^{-1}$. 
\end{proof}

\begin{remark}
We only need the case $m = 1$ in this paper. 
In this case, the lattice $M\subset V = E$ can be seen as an ideal of $E$, 
and we can work without translating into the language of hermitian lattices, 
although we do not go into the details.  
\end{remark}

\subsubsection{Overlattices}
We then explain how to list all lattices with similitude $(L, t)$ 
such that $(M_L, t)$ is isometric to a prescribed one.
An \emph{overlattice} of an integral lattice with similitude $(M, \hat{t})$ is a 
pair $(L, \hat{t})$ consisting of an overlattice $L\subset M^\vee$ of $M$ and the similitude 
obtained as the linear extension of $\hat{t}$, which is again denoted by $\hat{t}$, 
such that $\hat{t}(L) \subset L$. 
Every integral $(L, t)\in \cF_{S^m}$ is an overlattice of $(M_L, t)$. 
This implies that there are finitely many lattices with similitude $(L, t)$ 
such that $(M_L, t)$ is isometric to a prescribed one. 
In fact, the (integral) overlattices of a lattice $M$
are in bijection with the isotropic subgroups of the discriminant group $D_M$, 
and $D_M$ is a finite group (see e.g. \cite[\S 1.1.3]{Kondo_K3book}). 
Using this correspondence, it is possible to compute all overlattices of a 
given lattice (with similitude). 
In practice, we can narrow down the scope of isotropic subgroups using the conductor of $\o$. 

\begin{proposition}\label{prop:sandwich_of_L}
Let $(L, t)\in \cF_{S^m}$ be an integral lattice with similitude. 
Then we have the following sandwich of $L$ by $\cO_E$-modules:
\[ M_L \subset  L \subset M_L^\vee \cap \frakf^{-1}M_L. \]
\end{proposition}
\begin{proof}
Since $L$ is an overlattice of $M_L$, we have $M_L \subset L \subset M_L^\vee$. 
The inclusion $L \subset \frakf^{-1}M_L$ follows by \eqref{eq:fL<M_L}.  
\end{proof}

As a result, overlattices of a given $(M, s)$ that we need to compute 
correspond to isotropic subgroups of $(M^\vee \cap \frakf^{-1}M)/M \subset  D_M$. 
In \cite[\S 5.2]{Fieker-Hofmann-Sircana2019} Fieker, Hofmann and Sircana show how to efficiently enumerate subgroups of a finite abelian group $A$ stable under an isomorphism $a$. We modify their algorithm to enumerate isotropic subgroups only. This is necessary because the number of stable subgroups was impractically large in some cases. 

We remark that not every overlattice $(L, s)$ of a given $(M, s)\in \cF_{S^m}$ 
satisfies $(M_L, s) \cong (M, s)$ since $(L, s)$ may contain 
an $\cO_E$-module properly containing $M$ (even if it corresponds to 
an isotropic subgroup of $(M^\vee \cap \frakf^{-1}M)/M$). 
Furthermore, two different overlattices can be isometric.  
The following proposition summarizes this situation. 

\begin{proposition}
Let $(M,s)\in \cF_{S^m}$ be an integral lattice with similitude such that $M$ is 
an $\cO_E$-module. Let $\cE_M$ denote the set of overlattices 
$(L,s)$ of $(M,s)$ such that $L$ contains no $\cO_E$-submodule properly containing $M$. 
For each integral $(L,t)\in \cF_{S^m}$ and each isomorphism 
$\phi\colon (M_L, t) \to (M,s)$, we have an overlattice
$(\phi(L), s)$ of $(M,s)$ (after extending $\phi$ over $\bQ$). 
This correspondence gives rise to a bijection  
\[\{ (L, t)\in \cF_{S^m} 
\mid \text{$L$ is integral and $(M_L, t) \cong (M,s)$} \}/_{\cong} 
\xleftrightarrow{1:1} \cE_M/O(M,s).\] 
\end{proposition}
\begin{proof}
Straightforward.
\end{proof}

\begin{remark}
Note that $O(M,s) = U(M,h)$. 
In order to obtain an effective algorithm out of this proposition, we need to be able to compute the image of the natural map
\[U(M,h) \to \text{GL}_\o (\mathfrak{f}^{-1} M / M).\]
If $(M,h)$ is definite, the image can be derived by first calculating $U(M,h)$ using an algorithm of Plesken--Souvignier \cite{plesken-souvignier}. If $M$ is indefinite, one may use the strong approximation theorem via hermitian Miranda--Morrison theory \cite{brandhorst-hofmann}. 
In the cases relevant for us, $M$ is of rank one, and therefore $U(M,h)$ is the kernel of the norm map to the fixed field $K$ for the unit groups
\[\O_E^\times \to \O_K^\times, \quad \gamma \mapsto \gamma \sigma(\gamma),\]
which is readily computable. 
\end{remark}

Putting everything together, we obtain an algorithm for enumerating 
integral lattices with similitude in $\cF_{S^m}$. 

\begin{algorithm}\label{alg:enum_F_S^m}
Take as input an ideal $I$ of $\o$. 
The output is a complete set of representatives of the isomorphism classes of 
lattices with similitude $(L,t)\in \cF_{S^m}$ with  
\[I L^\vee \subset L\subset L^\vee.\]
\begin{steps}
\item Form the finite list $\mathcal{M}$ of hermitian $\O_E$-lattices $M$ such that 
\[\mathfrak{f} I M^\vee \subset M \subset M^\vee.\]
\item For each lattice $M \in \mathcal{M}$, 
compute a complete set of representatives of the $O(M,s)$-orbits 
($s$ is the multiplication by $\alpha$) of 
\[\{  L \mid \text{$L$ is a $\o$-module, 
$M \subset L \subset \mathfrak{f}^{-1} M \cap M^\vee$, $M_L = M$, and $L$ is integral} \}\]
and add the ones with $I L^\vee \subset L$ to the list of output lattices.
\end{steps}
\end{algorithm}
\subsection{The case of the transcendental parts}
In this subsection, we assume that $S(X)$ is a $q$-fractional Salem polynomial, 
and focus on the case where $(L, \hat{t})\in \cF_S$ is a
Frobenius lattice.
We use the same notation as in the previous subsection: 
$E = \bQ[X]/(S)$, $\alpha = X + (S)$, and so on. 
Let $\lambda\in \bQ_p$ be the unique root of $S(X)$ with $v_p(\lambda) > 0$. 
Then $S(X)$ decomposes over $\bQ_p$ as $S(X) = S_0(X)S_1(X)S_2(X)$, where
$S_0(X) = X- \lambda^{-1}$, $S_2(X) = X- \lambda$, and $S_1(X)$ is the remaining factor. 

We need a precise characterization of when $(L, \hat{t})\in \cF_S$ is a Frobenius lattice, 
i.e., satisfies \eqref{M4}.    
The following lemma gives us a practical criterion to test for it. 

\begin{lemma}\label{lem:is_frobenius_lattice}
Let $(L, t)\in \cF_S$ be an integral lattice with similitude. 
Then the following conditions are equivalent:  
\begin{enumerate}
\item $(L, t)$ is a Frobenius lattice; 
\item $L/\hat{t}(L) \cong (\bZ/q\bZ)^{\rk L - 2}\oplus \bZ/q^2\bZ$;
\item the elementary divisors of a matrix representing $\hat{t}$ 
in a given $\bZ$-basis of $L$ are $(1,q,\dots,q,q^2)$.
\end{enumerate}
\end{lemma}
\begin{proof}
For a prime $\ell$, we write $L_\ell = L\otimes \bZ_\ell$. 

(i) $\implies$ (ii). Since $\hat{t}(L_\ell) = L_\ell$
for every prime $\ell \neq p$, we have $L/\hat{t}(L) \cong L_p/\hat{t}(L_p)$. 
By \eqref{M4}, it follows that 
\[
\begin{split}
 L_p/\hat{t}(L_p) 
& = (L^0 \oplus L^1 \oplus L^2)/
\left( L^0 \oplus q L^1 \oplus q^2 L^2\right) \\
&= \left(L^1/qL^1\right) \oplus \left(L^2 / q^2 L^2\right) \\
&\cong (\bZ/q\bZ)^{\rk L - 2}\oplus \bZ/q^2\bZ. 
\end{split}
\]

(ii) $\iff$ (iii) is well known.

(iii) $\implies$ (i). Put 
$L^i = \{ x\in L_p \mid S_i(t)x = 0 \}$ for $i = 0,1,2$. 
We have $\hat{t} L^0 = q\lambda^{-1} L^0 = L^0$ and $\hat{t} L^2 = q\lambda L^2 = q^2 L^2$.
Note that the characteristic polynomial $\hat{S}(X)\in \bZ_p[X]$ of $\hat{t}$ decomposes
as $\hat{S}(X) = (X - q\lambda^{-1})\hat{Q}(X)$ for some $\hat{Q}(X)\in \bZ_p[X]$. 
Hence $\pi_0\coloneqq\hat{Q}(q\lambda^{-1})^{-1}\hat{Q}(\hat{t}):L_p\to L_p$ defines the projection to
$L^0$. Since $q^2\hat{t}^{-1}(L_p)\subset L_p$ by the assumption,
a similar argument gives the projection $\pi_2:L_p\to L_p$ to $L^2$.
Hence
\[ L_p = \im(\pi_0 + \pi_2)\oplus \im(\id - (\pi_0 + \pi_2)) 
= L^0\oplus L^2 \oplus L^1. \]
Finally, since $\hat{t} L^0 = L^0$ and $\hat{t} L^2 = q^2L^2$,
the assumption shows that the elementary divisors of $\hat{t}|_{L^1}$ are 
$(q,\ldots, q)$. Thus, $qL^1 = \hat{t}L^1$. 
Therefore, $(L, t)$ satisfies \eqref{M4}.  
\end{proof}

Our next aim is to refine \cref{prop:sandwich_of_L,prop:discriminant-annihilator} 
to fit the case of the transcendental parts of Frobenius lattices of K3 type.  
This refinement reduces the computational cost. 
To this end, we introduce a replacement for the conductor. 

Let $w_0$ and $w_2$ be the finite places of $E = \bQ[X]/(S)$ 
corresponding to the irreducible factors 
$S_0$ and $S_2$ of $S$ over $\bQ_p$, respectively. 
We define a subset $\frakf_\alpha \subset E$ by
\[\frakf_\alpha \coloneq \left\{ a\in E \mid a\cO_E \subset \sum_{i=0}^{2d-1}\bZ \alpha^{i} \right\}, \]
where $2d$ is the degree of $S(X)$. 
Then $\frakf_\alpha$ is the largest fractional ideal of $E$ contained in 
the $\bZ$-submodule $\sum_{i=0}^{2d-1}\bZ \alpha^{i}$. 
In particular, we have $\frakf_\alpha \subset \bZ[\alpha]$. 

\begin{remark}
We cannot define the fractional ideal $\frakf_\alpha$ as the subset 
\[\left\{ a\in E \mid a\cO_E \subset \bZ[\alpha] \right\} \]
in the same way as the usual conductor because $\alpha$ is not integral. 
In fact, this subset is the largest $\cO_E$-module contained in $\bZ[\alpha]$, but 
not finitely generated over $\cO_E$. 
%
\end{remark}

\begin{definition} \label{def:frobenius-conductor}
The \emph{Frobenius conductor} is the ideal
\[ \frakg\coloneq \frakg_\alpha \coloneq 
\frakf_{\alpha} \frakp_{w_0}^{-w_0(\frakf_{\alpha})} \frakp_{w_2}^{-w_2(\frakf_{\alpha})}, 
\]
where we use the same symbol for a finite place and its corresponding normalized valuation. 
\end{definition}

Recall that there exists a natural isomorphism $E\otimes_\bQ\bQ_p \to \prod_{w\mid p} E_w$ 
determined by $a\otimes b \mapsto (\iota_w(a)b)_w$, 
where $w$ ranges over all places of $E$ above $p$, and $\iota_w:E\to E_w$ denotes the natural inclusion. 
We identify both sides via this isomorphism. 
Then $\bZ_p \cO_E = \prod_{w\mid p} \cO_{E_w}$, and hence, for any ideal $\fraka$ of $\cO_E$, 
we have 
\begin{equation}\label{eq:decomp-of-ideal}
\bZ_p \fraka 
= \prod_{w\mid p} \cO_{E_w} \iota_w(\fraka) 
= \prod_{w\mid p} \bZ_p \iota_w(\fraka).  
\end{equation}

Let $(L,\hat{t})\in \cF_S$ be a Frobenius lattice. 
For a prime $\ell$, we write $L_\ell = L \otimes\bZ_\ell$. 
Note that $L^i = \{ x\in L_p \mid S_i(\alpha).x = 0 \}$ for $i = 0,1,2$, 
and we have $\alpha L^1 = L^1$. 
By \eqref{eq:decomp-of-ideal}, it follows that 
\[ \fraka L_p 
= \fraka (L^0\oplus L^2\oplus L^1)
= \iota_{w_0}(\fraka)L^0\oplus \iota_{w_2}(\fraka)L^2\oplus \iota_{1}(\fraka)L^1 \]
for any ideal $\fraka$ of $\cO_E$, 
where $\iota_1:E\to \prod_{w\mid p, w\neq w_0, w_2} E_w$ is the natural map. 
For simplicity, we sometimes omit the symbols $\iota_w$ or $\iota_1$. 

Put $V = L\otimes \bQ$. For each prime $\ell$, we denote by 
$\iota_\ell:V\to V_{\bQ_\ell} \coloneq V\otimes \bQ_\ell$ 
the natural inclusion. We have 
\begin{equation}\label{eq:reconstruction_from_completions}
L = \{ x\in V \mid \iota_\ell(x) \in L_\ell \textup{ for every prime $\ell$} \}.
\end{equation}
This is a standard fact from module theory and 
holds for any finitely generated free $\bZ$-module in $V$.

Using the Frobenius conductor $\frakg$, we have a counterpart of \eqref{eq:fL<M_L}
as follows. Let $\ell$ be a prime. We have 
$(\frakg L)_\ell 
= \bZ_\ell \iota_\ell(\frakg L)
= \frakg L_\ell$
in $V_{\bQ_\ell}$. 
If $\ell \neq p$, then the inclusion $\frakf_\alpha \subset \ZZ[\alpha]$ implies that
\[ (\frakg L)_\ell 
= \frakg L_\ell
= \frakf_\alpha L_\ell
\subset \bZ[\alpha] L_\ell
= L_\ell.
\]
If $\ell = p$, we have 
\[ (\frakg L)_p 
= \frakg L_p
= \frakg (L^0 \oplus L^2 \oplus L^1)
= \cO_{E_{w_0}} L^0 \oplus \cO_{E_{w_2}} L^2 \oplus \frakf_\alpha L^1.
\]
Since both $\cO_{E_{w_0}}$ and $\cO_{E_{w_2}}$ are equal to $\bZ_p$, 
and $\frakf_\alpha \subset \bZ[\alpha]$, we obtain
\[ (\frakg L)_p 
\subset \bZ_p L^0 \oplus \bZ_p L^2 \oplus \bZ[\alpha] L^1
= L^0 \oplus L^2 \oplus L^1
= L_p.
\]
Hence, it follows from \eqref{eq:reconstruction_from_completions} that 
$\frakg L \subset L$. By the maximality of $M_L$, we obtain
\begin{equation}\label{eq:gL<M_L}
\frakg L \subset M_L.
\end{equation}

\begin{proposition}[cf. {\cref{prop:discriminant-annihilator}}]
\label{prop:discriminant-annihilator:Frob}
Let $(L,t)\in \cF_{S}$ be a Frobenius lattice,  
and let $I$ be a discriminant annihilator for $(L,t)$. 
Then $\frakg I$ is a discriminant annihilator for $(M_L,t)$. 
\end{proposition}
\begin{proof}
By an argument similar to that in \cref{prop:discriminant-annihilator}, together with
\eqref{eq:gL<M_L}, it follows that
\[ \frakg I M_L^\vee \subset \frakg (I\cO_E L^\vee) 
\subset \frakg(\cO_E L) \subset \frakg(\frakg^{-1}M_L) = M_L. \]
\end{proof}

\begin{proposition}[cf. {\cref{prop:sandwich_of_L}}]
\label{prop:sandwich_of_L:Frob}
Let $(L,t)\in \cF_{S}$ be a Frobenius lattice. 
We have the following sandwich of $L$:
\[ M_L \subset  L \subset M_L^\vee \cap \frakg^{-1}M_L. \]
\end{proposition}
\begin{proof}
Since $L$ is an overlattice of $M_L$, we have $M_L \subset L \subset M_L^\vee$. 
The inclusion $L \subset \frakg^{-1}M_L$ follows by \eqref{eq:gL<M_L}. 
\end{proof}

By \cref{prop:sandwich_of_L:Frob,prop:discriminant-annihilator:Frob}, 
in \cref{alg:enum_F_S^m}, we can replace the conductor $\frakf$ by 
the Frobenius conductor $\frakg$ to collect Frobenius lattices $(L, t)$. 
We close this section by giving an example of a discriminant annihilator
for the transcendental part of a Frobenius lattice of K3 type.

\begin{example}\label{ex:disc-ann}
Let $(\Lambda, F)\in \cL$ be a Frobenius lattice of K3 type. Put 
$T = T(\Lambda, F)$ and $\NS = \NS(\Lambda, F)$. 
Let $\overline{C}(X)\in \bZ[X]$ be the minimal polynomial of 
the isometry $q^{-1}F|_{\NS}$ of $\NS$, and put $c = \overline{C}(\alpha)$. 
Then the ideal 
$I = (\o c)(\frakp_{w_0}^{-w_0(c)} \cap \o)(\frakp_{w_2}^{-w_2(c)} \cap \o)$ 
of $\o$ is a discriminant annihilator for $T$. 
In fact, since the minimal polynomial annihilates the action of $q^{-1}F|_{\NS}$ on $D_{\NS}$, 
the compatibility \eqref{eq:compatibility} yields 
$\overline{C}(q^{-1}F|_{T})D_{T} = 0$, which is equivalent to $c D_T = 0$.  
For every prime $\ell \neq p$, both $\frakp_{w_0}$ and $\frakp_{w_2}$ 
are units after localization, so $I (D_T)_\ell =c (D_T)_\ell = 0$. 
Equivalently, $I T_\ell^\vee \subset T_\ell$. 
At $p$, since $(D_T)_p = D_{T^1}$, the relation $c(D_T)_p = 0$ means that 
$c (T^1)^\vee \subset T^1$. Hence 
\[ I T_p^\vee = I (T^0\oplus T^2)^\vee \oplus I (T^1)^\vee 
= I(T^0\oplus T^2) \oplus c (T^1)^\vee 
\subset  (T^0\oplus T^2) \oplus T^1 = T_p. 
\]
Therefore, $I T^\vee\subset T$ globally.

Put $I'\coloneq (\cO_E c)\frakp_{w_0}^{-w_0(c)}\frakp_{w_2}^{-w_2(c)}$. 
A similar local argument gives $I'T^\vee\subset\cO_E T$. 
Then, the same proof as for \cref{prop:discriminant-annihilator:Frob} 
shows that $\frakg I'$ is a discriminant annihilator for $M_T$, 
which might properly contain $\frakg I$. 
\end{example}

\section{The N\'eron--Severi part}\label{sec:Neron--Severi-part}
In this section, we deal with N\'eron--Severi parts and prove \cref{th:has-rational-pts}.

\subsection{Algorithmic classification of lattices with isometry}
Recall that two lattices with isometry $(A,a)$ and $(B,b)$ are \emph{isomorphic} 
if there is an isometry $f \colon A \to B$ with $f \circ a= b \circ f$. 
They are said to be \emph{locally isomorphic} if  
$(A \otimes \bZ_\ell, a \otimes \id_{\bZ_\ell}) \cong (B \otimes \bZ_\ell, b \otimes \id_{\bZ_\ell})$ 
for all primes $\ell$ and 
$(A \otimes \bR, a \otimes \id_{\bR}) \cong (B \otimes \bR, b \otimes \id_{\bR})$.

In \cite{brandhorst-hofmann} the first author and Hofmann describe algorithms to compute a complete set of representatives of the isomorphism classes of lattices with isometry of finite order with given minimal and characteristic polynomial.

For the reader's convenience we sketch the approach of \cite{brandhorst-hofmann}.
Let $\mu \in \ZZ[X]$ be a squarefree and $*$-symmetric polynomial. 
Then 
\[E\coloneq\bQ[X,X^{-1}]/(\mu(X)) = \bQ[X]/(\mu(X))\]
is an \'etale $\bQ$-algebra and $X \mapsto X^{-1}$ defines an involution $\sigma$ on $E$. Let $\mu = \prod_{i \in I} \mu_i$ be the prime factorization of $\mu$. 
Set $E_i = \bQ[X]/(\mu_i)$. For the decomposition $E \cong \prod_i E_i$ we let $e_i$ be the corresponding system of primitive idempotents. 
Set $\Delta = \ZZ[X]/(\mu)$ and $\Gamma = \prod_{i \in I} \Delta e_i\cong \prod_{i \in I}\bZ[X]/(\mu_i)$
and let 
$\frakf_\Gamma=\{a \in \Delta \mid a \Gamma \subseteq \Delta\}$ denote the conductor of $\Gamma$ in $\Delta$. 

Let $L$ now be a lattice and $f\in O(L)$ with minimal polynomial $\mu$. Set 
\[L_i  = L \cap e_i L \quad \mbox{ and } \quad L_i' = e_i \Gamma L.\]
Then we have a sandwich
\[\frakf_\Gamma L \subseteq \bigperp_{i \in I} L_i \subset L \subset \Gamma L = \bigperp_{i \in I} L_i'\]
where the two outermost pieces are direct sums of hermitian lattices over orders in number fields.
In our situation, $f$ is of finite order and thus $\mu \mid X^n-1$. Then the orders $\Delta e_i$ are maximal and the corresponding hermitian lattices are well understood and amenable to computations.

\begin{example}
Let $L= \ZZ^2$ and $f \in O(L)$ defined by $f((x,y)) = (y,x)$. Then $\mu=X^2-1$, 
$e_1=(X+1)/2$, $e_2 = -(X-1)/2$, $\Delta = \bZ[X]/(X^2-1)$, $\Gamma \cong \bZ[X]/(X-1) \times \bZ[X]/(X+1) \cong \bZ \times \bZ$, $L_1 = (1,1)\ZZ$, $L_2=(1,-1)\ZZ$, $\frakf_\Gamma = 2\Gamma$, $\Gamma L =\{(a/2,b/2) \mid a \equiv b \mod 2\}= \frac{1}{2} (L_1 \perp L_2)$. 
\end{example}

Via the conductor one obtains bounds on the discriminants of the $L_i$.
When one wants to enumerate (possibly local) isomorphism classes of lattices with isometry and fixed minimal polynomial, one first lists all possibilities for the $L_i$ and then takes (suitable) overlattices of $\frakf_\Gamma L$ modulo the action of the unitary groups 
$\prod_{i\in I} U(L_i)$. In practice, one successively takes equivariant primitive extensions.

\begin{remark}
The computations done for this work are several orders of magnitude larger than the ones for \cite{brandhorst-hofmann}. The original algorithms have no chance to terminate within a reasonable time. We therefore improved them for this work:
\begin{itemize}
    \item we view the discriminant group as a $\ZZ[X]$-module and give its annihilator ideal as input to the algorithm in order to discard candidates early on;
    \item we first enumerate the local isomorphism classes and then split each local one into global ones. Both steps are computationally very challenging; 
    \item the splitting into global isometry classes requires the enumeration of genera of definite integer lattices via Kneser's neighbor method. For this the first author and Hofmann improved the canonization and automorphism group computation of definite integer lattices in \cite{Brandhorst-Hofmann2026}. 
\end{itemize}
The improvements are available in OSCAR \cite{OSCAR}. 
\end{remark}

Now we turn to effectivity:
\begin{lemma}
Let $(\Lambda, F)$ be a Frobenius lattice of K3 type. 
Put $\NS = \NS(\Lambda, F)$ and $f = q^{-1}F$.  
Let $\NS_f\coloneq\left(\NS^f\right)^\perp \subset \NS(\Lambda, F)$ 
denote the coinvariant lattice of $f|_{\NS(\Lambda,F)}$. 
Then $(\Lambda, F)$ is effective if and only if 
$\NS_f$ is negative definite and does not contain any roots,
where a root is an element with self-inner product $-2$.
\end{lemma}
\begin{proof}
Note that $f|_{\NS}$ has finite order by \cref{prop:NS}.
Suppose that $(\Lambda, F)$ is effective. 
Averaging an interior vector of a Weyl chamber preserved by $f$ over the cyclic group generated by $f|_{\NS}$ gives an invariant interior vector $v$.
Since $v^2 >0$ and $\NS$ has signature $(1, \rk \NS -1)$, the coinvariant lattice
$\NS_f$ is negative definite. Furthermore, no root can lie in $\NS_f$ 
because $v$ would then lie on a root hyperplane.  

Conversely, suppose that $\NS_f$ is negative definite and does not contain any roots. 
Then $\NS^f$ contains a vector $v$ with $v^2 > 0$. 
Since $\NS_f$ contains no root, for any root $r\in \NS$, 
the intersection $(\NS_\bQ^f\cap r^\perp)$ is a proper subspace of $\NS_\bQ^f$.   
Since the root hyperplanes are locally finite in the positive cone,
the vector $v$ can be chosen outside all root hyperplanes. 
Then the Weyl chamber containing $v$ is preserved by $f|_{\NS}$. 
\end{proof}

The \emph{center density} of a definite lattice $L$ is defined as 
$\delta(L) = \rho(L)^{\rk L}/ \sqrt{|\det L|}$ 
where $\rho(L)= 2^{-1}\sqrt{\min(L)}$
and $\min(L)= \min \{ |x^2| \mid x \in L, x\neq 0 \}$.
There are upper bounds $b_n$ on the center density of a lattice depending only on its rank $n$. We use the list of Henry Cohn \cite{cohn-bounds}.

\begin{proposition}[Sphere packing test]
If $L$ is a negative definite even lattice of rank $n$ with 
$|\det(L)|^{-1/2} > b_n$, then $L$ represents $-2$.
\end{proposition}
\begin{proof}
If $L$ does not represent $-2$, then 
$\min(L) \geq 4$, hence 
$\rho(L)\geq 1$ and 
\[b_n \geq \delta(L) \geq |\det(L)|^{-1/2}.\]
The claim follows by contraposition.
\end{proof}

\subsection{Algorithmic classification of effective Frobenius lattices of K3 type}
We now prove \cref{th:has-rational-pts}. 
Recall that for a K3 surface $X$ over $\bF_q$ and its corresponding triple 
$(\Lambda, F, \cK)$, the number of rational points is given by the 
formula \eqref{eq:point-count}. 

\begin{algorithm}\label{alg:non-positive-point-count}
Fix a prime power $q$. 
The output of the following algorithm 
is a complete set of representatives of the isomorphism classes of 
effective Frobenius lattices of K3 type $(\Lambda, F)$ such that 
$1 + \Tr(F) + q^2 \leq 0$. 
\begin{steps}
\item Enumerate all $q$-fractional Salem polynomials of degree $2 \leq 2d \leq 20$ 
satisfying the trace bound \eqref{eq:trace-bound}. 
\item Fix a $q$-fractional Salem polynomial $S(X)$ from Step 1. 
Compile the list consisting of all products of cyclotomic polynomials $C(X)$ 
such that $1 + q\Tr S(X)C(X) + q^2 \leq 0$, $C(1) = 0$, and $qS(-1)C(-1)$ is a square. 
\item For each such $C(X)$, compile the finite list of all candidates for 
transcendental parts $(L, \hat{t})$ such that $q^{-1}\hat{t}$ has characteristic polynomial 
$S(X)$ as explained in \S\ref{sec:transcendental-part}.  
\item For each such $(L,\hat t)$ and $C(X)$, compute the list of all local isometry classes 
of lattices with isometry $(N,n)$ compatible with $(L, \hat t)$, 
and apply the sphere packing test to exclude some which cannot be effective. 
\item For each $(N,n)$ compute genus representatives of the negative definite lattices $\ker(\id+n)$. Discard the ones representing $-2$. 
\item Compute a complete set of representatives of the global isometry classes of the lattices $(N, n)$
and check which ones are effective.
\item Compute the double coset \eqref{eq:double-coset}. 
Using each anti-isometry that represents an element of the double coset, 
glue $(N, n)$ and $(L, \hat{t})$ to obtain an effective Frobenius lattice of K3 type
$(\Lambda, F)$ such that $1 + \Tr(F) + q^2 \leq 0$.  
\end{steps}
\end{algorithm}

\begin{proof}[Proof of \cref{th:has-rational-pts}]
Let $q$ be a power of a prime.
If $q>19$ then it follows from the Weil conjectures that every K3 surface 
over $\bF_q$ has an $\bF_q$-point. 
Suppose that $13 \leq q \leq 19$. We carry out \cref{alg:non-positive-point-count}. 
Note that \cref{th:criteria-for-having-pts}, which we prove 
in \S\ref{sec:geometric-approach} independently, shows that 
a K3 surface over $\bF_q$ has an $\bF_q$-point if the $q^{-1}F|_{\NS}$-invariant lattice of $\NS(\Lambda, F)$ represents $-2$ or $0$,
where $(\Lambda, F)$ is the corresponding lattice with similitude of K3 type.  
Since it is determined only by its genus whether a lattice represents $0$ or not, 
we can add the test whether the $n$-invariant part of $(N,n)$ represents $0$ 
at Step 4 in \cref{alg:non-positive-point-count}. Since the $n$-invariant part is indefinite (if it is of rank at least two), it is frequently unique in its genus, hence checking for representability of $(-2)$ is again a local condition. 
These tests reduce the computational cost. No candidate passed Step 5 for $q \in \{13,16,17,19\}$.
By \cref{th:Taelman2020}, this proves that there
are no pointless K3 surfaces over $\bF_q$. 
This completes the proof. 
\end{proof}

\begin{remark}
When an effective Frobenius lattice is induced by an ordinary K3 surface $X/\FF_q$, its point counts $a_m=\# X(\bF_{q^m})$ are nonnegative and satisfy $a_{mn}\geq a_m$ for all $m,n\geq 1$.

This is certainly the case if \condAst holds, e.g. when $p\geq 5$ and $\rho(X)\geq 12$. One may also assume it if one is just interested in disproving the existence of pointless K3 surfaces. 
Since we are also interested in testing the conjecture on potential semi-stable reduction \condAst, we did not assume this. 
\end{remark}

\begin{remark}
The computations took (very roughly) 100 years of CPU time on a cluster. 
The results of the genus enumerations are available as part of the database \cite{Muller-Brandhorst2026:db}.
The implementation and the computational data are available at \cite{ancillary-files}.

After Step 4, for $q=19,17,16,13$, the number of candidates remaining is $15$, $29920$, $96$, $\geq 1682165$. For $q=13$ we only stored the $\NS$-parts if there are more than $10000$ transcendental parts with the same $S(X)$. Most the Frobenius lattices of K3 type for $q=16$ cannot come from a K3 surface because their virtual point counts are negative.  (If \condAst holds, then we know that they are not effective.)
\end{remark}


We illustrate the case $q=19$. 

\subsubsection*{The characteristic polynomial}
By \Cref{table:q-salem-small-trace}, there are 38 $19$-fractional Salem polynomials
of degree at most $20$ satisfying \eqref{eq:trace-bound}.
They are all of degree at most $4$. Since potential semistable reduction \condAst
is known for K3 surfaces over $\mathbb{F}_{19}$ and geometric Picard number $\geq 18$, 
we will not find effective Frobenius lattices of K3 type with negative point count
(see e.g. \cite[Theorem 4.6]{Taelman20}). 
Therefore, it suffices to search for cyclotomic polynomials $C(X)$ such that 
$1 + q\Tr S(X)C(X) + q^2 =0$. 
There are only five polynomials 
$S(X)C(X)$ satisfying this condition, namely, $S(X)C(X)$
where $C(X) = (X-1)(X+1)^{21-\deg S}$ and 
$S(X)$ is among the following 5 $19$-fractional Salem polynomials: 
\begin{eqnarray*}
S_1(X) = X^{2} + \frac{20}{19} X + 1\\
S_2(X) = X^{4} + \frac{58}{19} X^{3} + \frac{79}{19} X^{2} + \frac{58}{19} X + 1\\
S_3(X) = X^{4} + \frac{58}{19} X^{3} + \frac{80}{19} X^{2} + \frac{58}{19} X + 1\\
S_4(X) = X^{4} + \frac{58}{19} X^{3} + \frac{81}{19} X^{2} + \frac{58}{19} X + 1\\
S_5(X) = X^{4} + \frac{58}{19} X^{3} + \frac{82}{19} X^{2} + \frac{58}{19} X + 1.
\end{eqnarray*}
In what follows, we use the same notation as in \S\ref{sec:transcendental-part}: 
$E = \bQ[X]/(S)$, $\alpha = X + (S)\in E$, $\hat{\alpha} = q\alpha\in \cO_E$, and so on. 
Note that the minimal polynomial of the N\'eron--Severi part will be 
$\overline{C}(X) = X^2 -1$. Putting 
$I' = (\alpha^2 - 1)\frakp_{w_0}^{-w_0(\alpha^2 - 1)}\frakp_{w_2}^{-w_2(\alpha^2 - 1)}$, 
we can take the ideal $\frakg I'$ as a discriminant annihilator for $M_T$, 
where $T$ is the transcendental part (\cref{ex:disc-ann}). 
We begin with the case $S = S_2$. 

\subsubsection*{The transcendental part}
The number field $E = \bQ[X]/(S_2)$ has (relative) class number $6$.
The order $\ZZ[\hat \alpha]$ is not maximal and its conductor is of norm $19^6$,  
but its Frobenius conductor $\mathfrak{g}$ is trivial. 
Therefore any Frobenius lattice in $\cF_S$ comes from a hermitian $\cO_E$-lattice.

The fractional ideal $(\alpha^2 - 1)\O_E$ factors as $\mathfrak{p}_{233} \mathfrak{p}_{w_0}^{-2}$, 
where $\mathfrak{p}_{233}$ is a prime of norm $233$. 
So we know that $L = M_L$ satisfies $\mathfrak{p}_{233}L^\vee \subset L$. 
By a calculation with OSCAR, there is, up to isometry, a unique genus $\mathfrak{t}$ 
consisting of hermitian lattices $T$ such that $\mathfrak{p}_{233}T^\vee \subset T$ 
and whose trace lattice has signature $(2,2)$. Since the relative class number of $E$ is $6$  
(and $\dim_E T\otimes \bQ =1$), the genus $\mathfrak{t}$ consists of $6$ isometry classes 
(cf. \cite[Proposition 2.9]{brandhorst-cattaneo}). Fix any $(T,h) \in \mathfrak{t}$.

A direct calculation shows that its trace lattice $(T,\Tr_{E/\bQ} \circ h)$ is in the genus $\II_{2,2} 233^1$ and that the isometry $t$ defined by multiplication with $\alpha$ acts as the identity on the discriminant group $D_T$.

\subsubsection*{The N\'eron--Severi part}
By \cite[Proposition 1.14.1]{nikulin1980_forms} (or a computation with OSCAR) $T$ admits a unique embedding into the K3 lattice 
$\Lambda$. Its orthogonal complement $\NS$ in $\Lambda$ satisfies $D_{\NS}\cong D_T(-1)$, 
is even and of signature $(1,17)$. Therefore it is in the genus $\II_{1, 17} 233^1$  (since $-1$ is a square mod $233$ and the exponent is odd, the sign of the exponent of $233$ is the same). 

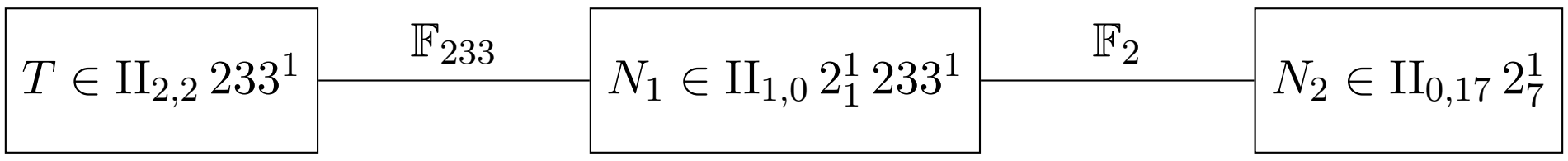
\begin{figure}[!htbp]
  \tikzstyle{block} = [draw, rectangle, minimum height=3em, minimum width=3em]
  \tikzstyle{virtual} = [coordinate]
  \centering
  \begin{tikzpicture}[auto]
    \node [block] at (-5,0)  (T)     {$T \in \II_{2,2}\, 233^1$};
    \node [block] at (0,0) (C1)     {$N_1\in \II_{1,0}\,2^1_1 \,233^1$};
    \node [block] at (5,0) (C2) {$N_2\in \II_{0,17}\,2^1_7$};
    \draw [-] (T) -- node {$\mathbb{F}_{233}$} (C1);
    \draw [-] (C1) -- node {$\mathbb{F}_{2}$} (C2);
  \end{tikzpicture}
  \caption{Glue diagram for $S_2$}\label{fig:glueS2}
 \end{figure}

Suppose there exists an isometry $n \in O^+(\NS)$ of order $2$ with characteristic polynomial 
$(X-1)(X+1)^{17}$ which acts trivially on the discriminant group $D_{\NS}\cong D_T(-1)$. 
Set $N_i \coloneq \ker c_i(n)$, where $c_1(X) = X-1$ and $c_2(X) = X+1$. 
We have $N_1 \oplus N_2 \subset \NS$ and know that $N_1$ is of signature $(1,0)$ and $N_2$ of signature $(0,17)$. 
Since the reduced resultant of $X-1$ and $X+1$ is $2$, we know that $2\NS \subset N_1 \oplus N_2$. Because $1 = m\coloneq\min \{\rk N_1,\rk N_2\}$, we know that $[\NS:N_1\oplus N_2] \mid 2^m=2^1$. 
We know that $-233 = \det \NS = \det N_1 \det N_2 / [\NS:N_1\oplus N_2]^2$ and therefore that $233$ divides $\det N_1$ or $\det N_2$ and in fact that $D_{\NS}$ embeds (equivariantly) into $D_{N_i}$ for $i=1$ or $i=2$.
But $n$ acts trivially on $D_{\NS} \cong \bF_{233}$ and so $i=1$. Because $\NS$ is even, so is $N_1$, which is only possible if $2 \mid \det N_1$. Hence $\det N_1=2\cdot 233$, $\det N_2=-2$ and $[\NS:N_1\oplus N_2]=2$. Thus $N_1 \cong [2 \cdot 233] \in \II_{1,0}2^1_1 233^1$. Then we know that $(D_{N_1})_2 \cong D_{N_2}(-1)$ and can conclude that $N_2 \in \II_{0,17} 2^1_7$. The situation is summarized in \cref{fig:glueS2}.
These $3$ lattices glue to a Frobenius lattice of K3 type, but any lattice $N_2 \in \II_{0,17} 2^1_7$ does not pass the sphere packing test, and therefore, the Frobenius lattice in question is not effective. \\

For $S_3$ the conductor of $\bZ[\hat \alpha]$ is of norm $3^2 19^6$ and the Frobenius conductor $\mathfrak{g}$ is of norm $3^2$. 
A computer calculation shows that there are $6$ genera of hermitian lattices $M$ with $\mathfrak{g}I' M^\vee \subset M$. Since the Frobenius conductor $\mathfrak{g}$ is non-trivial, we have to consider 
overlattices $T \supset M=M_T$.
A computer calculation following \Cref{alg:enum_F_S^m} 
(we use $\frakg$ instead of $\frakf$) shows that there are $12$ Frobenius lattices $T$ with $(\alpha^2-1)D_T=0$.
They lie in the following $8$ genera. 
\[\II_{2,2} 2^{\pm 2} 13^{\mp}, \quad \II_{2,2} 2^{\pm 2}3^2 13^{\mp},\quad  \II_{2,2}2^2 9^\pm 13^\mp,\quad  \II_{2,2}2^{-2} 9^\pm 13^\mp.\] 
As before one obtains $8$ possible genera for $\NS$.
In each case, $D_t|_{(D_T)_\ell} = \id_{(D_T)_\ell}$ for $\ell \neq 2$, 
and $t$ acts with order two on $(D_T)_2$. 
We know that $(D_T)_\ell(-1) \cong (D_{N_1})_\ell$ and that $\ell \nmid \det N_2$.
This excludes $\II_{2,2} 2^{\pm 2}3^2 13^{\mp}$ because as groups $(D_{T})_3 \cong \bF_3^2 \cong (D_{N_1})_3$ and the right hand side is too large to come from a lattice of rank $1$. 

From the reduced resultant of $x-1$ and $x+1$ we know that $[\NS:N_1\oplus N_2] \mid 2$ and hence $N_2$ has determinant dividing $-8$ and signature $(0,17)$. 
It is excluded by the sphere packing test. 

The cases $S_4$ and $S_5$ are excluded in a completely analogous way. 
For $S_4$ one obtains $|\det N_2|=6$, and for $S_5$ one obtains $|\det N_2| \in \{2,8\}$.

In the case $S_1$ a similar computation shows that $N_1\cong [116]$ and gives
$4$ possible genera for $N_2$, namely $\II_{0,19}4^{-1}_5 a$ where $a \in \{1, 3^{-2},9^{\pm 1}\}$. The first one is excluded by the sphere packing test. We enumerated the remaining $3$ genera with Kneser's neighbor method and obtained $709$, $1052$, and $1017$ distinct isometry classes, respectively (available in the database \cite{Muller-Brandhorst2026:db}). 
All of them represent $-2$. Therefore the corresponding Frobenius lattices of K3 type are not effective.

\section{Geometric approach}\label{sec:geometric-approach}
In this section, we prove \cref{th:criteria-for-having-pts}.



\subsection{Geometry of K3 surfaces over non-closed fields}
To prove \cref{th:criteria-for-having-pts}, we study the geometry of K3 surfaces over non-closed fields.
Let $k$ be a field. We fix an algebraic closure $\overline{k}$ of $k$, and 
let $k^\textrm{s}$ denote the separable closure of $k$ in $\overline{k}$. 
Set $G \coloneq \Gal(k^\textrm{s}/k)$. 
For a scheme $X$ over $k$, we write 
\[X^\mathrm{s} = X\times_k\Spec(k^\textrm{s}), \quad
\overline{X} = X\times_k\Spec(\overline{k}).
\]

We begin by recalling some facts from \cite{Brigh-Logan-vanLuijk2020}. 
Let $X$ be a K3 surface over $k$, and put $L \coloneq \Pic(X^\mathrm{s})$. 
The natural map $\Pic(X^\mathrm{s})\to \Pic(\overline{X})$ is an isomorphism; 
see \cite[Lemma 3.1]{Brigh-Logan-vanLuijk2020}. 
Note that the Galois group $G$ acts on $L$ isometrically in a natural way. 
We define the action of $G$ on $\rO(L)$ as 
$\sigma.t = \sigma\circ t \circ \sigma^{-1}$ for $\sigma\in G$ and $t\in \rO(L)$. 
Let $\Delta\coloneq \{\delta\in L\mid \delta^2 = -2 \}$, and  
let $W\subset \rO(L)$ be the subgroup generated by $\{ r_\delta \mid \delta\in \Delta \}$, 
where $r_\delta$ is the reflection orthogonal to $\delta$. 
The subgroup $W^G\subset W$ of $G$-invariant elements is denoted by $R_X$. 
We identify $\Pic(X)$ with its image under the natural injection
\[ \Pic(X) \to L^G = \Pic(X^\mathrm{s})^G. \]
Since the image is of finite index in $L^G$,
we have $\Pic(X)_\bR = L_\bR^G$. 
Let $\cC_X \subset \Pic(X)_\bR$ be the positive cone of $X$. 
Note that $R_X\subset W$ preserves $\cC_X$.

Let $\cI$ be the set of $G$-orbits $I$ of $(-2)$-curves on $X^\mathrm{s}$ of 
the following two types:
\begin{enumerate}
\item $I$ consists of disjoint $(-2)$-curves;
\item $I$ consists of disjoint pairs of $(-2)$-curves, each pair having 
intersection number $1$. 
\end{enumerate}
For a $G$-orbit $I = \{ E_1, \ldots, E_m \}$ of type (i), we define 
\[r_I \coloneq r_{E_1} \circ \cdots \circ r_{E_m}\in W.   \]
For a $G$-orbit $I = \{ E_1,E_1',  \ldots, E_m, E_m' \}$ of type (ii) with
$(E_i, E_i') = 1$, we define 
\[r_I \coloneq r_{E_1+E_1'} \circ \cdots \circ r_{E_m+E_m'}\in W.   \]
Note that $r_I$ does not depend on the way to label the members of $I$, and 
it belongs to $R_X$. 
It follows that $R_X$ is generated by $\{ r_I \mid I\in \cI \}$, and that 
for each $I\in \cI$, the restriction $r_I|_{L_\bR^G}:L_\bR^G\to L_\bR^G$ 
is the reflection orthogonal to 
\[ C_I \coloneq \sum_{E\in I} E \in L^G; \]
see \cite[Proposition 3.6]{Brigh-Logan-vanLuijk2020}. 
Let $H_I\subset L_\bR^G$ denote the hyperplane $C_I^\perp$ in $L_\bR^G$. 
Let $\Nef(X)\subset L_\bR^G$ be the nef cone of $X$, and define  
$\cD \coloneq \Nef(X)\cap \cC_X$. 

\begin{proposition}\label{prop:BLvL}
The cone $\cD$ is a fundamental domain for the action of $R_X$ 
on $\cC_X$, and this action is faithful.  
Moreover, every wall of $\cD$ is the hyperplane $H_I$ for some $I\in \cI$. 
\end{proposition}
\begin{proof}
See Proposition 3.7 and Remark 3.9 in \cite{Brigh-Logan-vanLuijk2020}. 
\end{proof}

The following proposition is the key to proving \cref{th:criteria-for-having-pts} (i),
which is itself remarkable. 

\begin{proposition}\label{prop:-2element_general}
Suppose that there exists $\delta\in \Pic(X)$ such that $\delta^2 = -2$. 
Then $X^\mathrm{s}$ contains either a $G$-invariant $(-2)$-curve or a $G$-invariant pair $\{E, E'\}$ of $(-2)$-curves such that $(E, E') = 1$ and $\sigma(E) = E'$ for some $\sigma\in G$. 
\end{proposition}
\begin{proof}
The reflection $r_\delta$ belongs to $R_X$ since 
$\sigma \circ r_\delta \circ \sigma^{-1} = r_{\sigma(\delta)} = r_\delta$. 
By \cref{prop:BLvL}, the action of $R_X$ on $L_\bR^G$ is faithful, and 
hence $r_\delta|_{L_\bR^G}$ is a reflection with respect to some hyperplane 
in $L_\bR^G$. Furthermore, since $\cD$ is the fundamental domain for the action of $R_X$ 
on $\cC_X$, the reflection $r_\delta|_{L_\bR^G}$ is conjugate in $R_X$ to 
a reflection with respect to some wall of $\cD$, which is $H_I$ for some $I\in \cI$. 
That is, $r_\delta|_{L_\bR^G}$ is conjugate to $r_I|_{L_\bR^G}$. 
This implies that $r_\delta$ is conjugate to $r_I$, again because the action 
of $R_X$ on $L_\bR^G$ is faithful.
Hence 
\[ \dim \ker(r_I + \id) =  \dim \ker(r_\delta + \id) = 1. \]
This dimension is the number of $(-2)$-curves in $I$ if $I$ is of type (i) and 
the number of pairs of $(-2)$-curves in $I$ if $I$ is of type (ii).  
This completes the proof. 
\end{proof}

We next consider the case where $\Pic(X)$ contains an isotropic element, 
that is, a nonzero element of self-intersection $0$.

\begin{lemma}[cf. {\cite[Remark 8-2.13]{Huybrechts_K3book}}]\label{lem:nef-line-bundle}
Suppose that $\Pic(X)$ contains an isotropic element. 
Then, there exists a nontrivial nef line bundle on $X$ of self-intersection zero. 
\end{lemma}
\begin{proof} 
Let $\delta\in \Pic(X)$ be an isotropic element. 
We may assume that $\delta$ is effective by replacing it by $-\delta$ if necessary. 
We fix an ample class $h\in \Pic(X)$. Then $(h, \delta)>0$. 
We show that if $\delta$ is not nef then there exists an isotropic element $\delta' \in \Pic(X)$ 
that is effective and satisfies $0 < (h, \delta') < (h, \delta)$. 
Suppose that $\delta$ is not nef. Since 
\[ \Nef(X) = \Nef(X^\mathrm{s}) \cap L_\bR^G 
= \{ x\in \overline{\cC_X} \mid (x, E) \geq 0 \text{ for any $(-2)$-curve $E$ on $X^\mathrm{s}$} \},  
\]
there exists a $(-2)$-curve $E_1$ on $X^\mathrm{s}$ such that $(\delta, E_1) < 0$. 
Let $I = \{ E_1, \ldots, E_m \}$ be the $G$-orbit of $E_1$, and put 
$C_I \coloneq \sum_{i=1}^m E_i$. Since the action of $G$ on $L$ is isometric, 
we have
\[ (\delta, C_I) = \sum_{i=1}^m (\delta, E_i) = m(\delta, E_1) < 0.  \]
We now claim that $C_I^2 < 0$. 
Suppose to the contrary that $C_I^2 \geq 0$. Then 
$C_I$ lies in the closure $\overline{\cC_X}$ of $\cC_X$ since it is effective. 
Note that $\delta$ also lies in $\overline{\cC_X}$. 
Put $\lambda \coloneq -(C_I^2 + 1)/(\delta, C_I) > 0$. 
Since $\overline{\cC_X}$ is a cone, we have $C_I + \lambda\delta \in \overline{\cC_X}$, 
and thus $0 \leq (C_I + \lambda\delta)^2$. 
On the other hand, 
\[ (C_I + \lambda\delta)^2 = C_I^2 + 2 \lambda(C_I, \delta) = -C_I^2 -2 < 0. \]
This is a contradiction. Hence $C_I^2 < 0$. 

For each $i$, set $d_i \coloneq (E_i, \sum_{j\neq i}E_j) \geq 0$. Then
$d_1, \ldots, d_m$ are the same, say $d$, since the action of $G$ on $L$ is isometric.
By the inequality
\[ 0 > C_I^2 = \sum_{i=1}^m(E_i^2 + d_i) = m(-2 + d),  \]
it follows that $d = 0$ or $1$.
If $d = 0$ then the $G$-orbit $I$ is of type (i), and
if $d = 1$ then the $G$-orbit $I$ is of type (ii).
Let $\delta' \coloneq r_I(\delta)$. Then $\delta'^2 = \delta^2 = 0$, and
$\delta'$ belongs to $\overline{\cC_X}$, which implies that $\delta'$ is effective,
and $(h, \delta')>0$. Moreover,
since $r_I|_{L_\bR^G}$ is the reflection orthogonal to $C_I$, we obtain
\[ (h, \delta')
= \left(h, \delta - \frac{2 (\delta, C_I)}{C_I^2} C_I \right)
= (h, \delta) - \frac{2 (\delta, C_I)}{C_I^2} (h, C_I) < (h, \delta)
\]
as required.
Since $(h, \Pic(X)) \subset \bZ$, repeating this process terminates and 
yields a nontrivial isotropic element that is nef. This completes the proof.
\end{proof}

\begin{proposition}\label{prop:genus-one-fibration}
Suppose that $\Pic(X)$ contains an isotropic element. 
Then $X$ admits a genus-one fibration $f:X\to B$ with 
$\overline{B} \cong \bP^1_{\overline{k}}$.     
\end{proposition}
\begin{proof}
By \cref{lem:nef-line-bundle}, 
there exists a nontrivial nef line bundle $L$ on $X$ with $(L,L) = 0$. 
Let $\overline{L}$ denote the line bundle on $\overline{X}$ obtained by
the base extension of $L$. 
Then, by \cite[Proposition 3.10]{Huybrechts_K3book}, the line bundle 
$\overline{L}$ is globally generated, and hence, so is $L$. 
Let $\varphi:X\to \bP^n_k$, where $n = \dim H^0(X, L)-1$, be the morphism defined by $L$, 
and let $Y$ be the image of $\varphi$. 
Then, the base extension $\overline{\varphi}:\overline{X}\to \bP^n_{\overline{k}}$ 
of $\varphi:X\to \bP^n_k$ to $\overline{k}$ is the morphism defined by $\overline{L}$, 
and its image is $\overline{Y}$. 
As in the proof of \cite[Proposition 3.10]{Huybrechts_K3book}, 
we have $\dim\overline{Y} = 1$, and the Stein factorization of $\overline{\varphi}$ 
yields a genus-one fibration $\overline{f}: \overline{X}\to \overline{B} \cong \bP^1$. 
Since the Stein factorization commutes with flat base extension 
(cf. \cite[Part 2, Section 37.53]{StacksProject}), 
that of $\varphi$ yields the desired genus-one fibration. 
\end{proof}

We will need the following lemma. 

\begin{lemma}\label{lem:quasiell_implies_ss}
Suppose that $X$ admits a quasi-elliptic fibration. 
Then $X$ is supersingular. 
In particular, if $k$ is a finite field, then $X(k)\neq \emptyset$. 
\end{lemma}
\begin{proof}
Let $X\to B$ be a quasi-elliptic fibration. Then, its base extension to 
$\overline{k}$ is also a quasi-elliptic fibration. 
This implies that $\overline{X}$ is uniruled, 
and thus, $\overline{X}$ is supersingular;
see e.g. \cite[\S2]{Kawabe2024} for the proof.
Hence $X$ is supersingular and not ordinary. 
In particular, if $k$ is a finite field, then $X(k)\neq \emptyset$ 
because $X$ is ordinary if and only if $\#X(k)\not\equiv 1 \bmod \mathrm{char}(k)$. 
\end{proof}

\subsection{Proofs}
Now, we prove \cref{th:criteria-for-having-pts}. 
Throughout this subsection, $k$ is a finite field, and $X$ is a K3 surface over $k$. 
Note that $k^\textrm{s} = \overline{k}$. 
We also remark that there is no nontrivial Severi--Brauer curve over $k$, 
that is, for any curve $C$ over $k$ with $\overline{C} \cong \bP_{\overline{k}}^1$
we have $C \cong \bP_{k}^1$. 
This follows from the fact that $\mathrm{Br}(k)=0$; see e.g. \cite[\S 4.5.1]{Poonen_rational-points}. 

\begin{proof}[Proof of \cref{th:criteria-for-having-pts} (i)]
Suppose that there exists $\delta\in \Pic(X)$ such that $\delta^2 = -2$. 
By \cref{prop:-2element_general}, 
$\overline{X}$ contains either a $G$-invariant $(-2)$-curve or a $G$-invariant pair $\{E, E'\}$
of $(-2)$-curves such that $(E, E') = 1$. 
Suppose first that $\overline{X}$ contains a $G$-invariant $(-2)$-curve $E$. 
Then, $E$ is defined over $k$ since it is $G$-invariant. 
Namely, $X$ contains a curve $C$ such that $C_{\overline{k}} = E$. 
Since there is no nontrivial Severi--Brauer curve over $k$, 
we obtain $C \cong \bP^1_{k}$. In particular, 
$C$ has a $k$-point, which belongs to $X(k)$. Hence $X(k)$ is non-empty. 

Suppose then that $\overline{X}$ contains a $G$-invariant pair $\{E, E'\}$
of $(-2)$-curves such that $(E, E') = 1$.  
Let $p\in X(\overline{k})$ be the unique intersection point of $E$ and $E'$. 
Then $p$ is invariant under $G$. 
Thus $p\in X(k)$, and $X(k)$ is non-empty. The proof is complete.  
\end{proof}

\begin{proof}[Proof of \cref{th:criteria-for-having-pts} (ii)]
Suppose that there exists a nonzero element $\delta\in \Pic(X)$ with $\delta^2 = 0$. 
Then $X$ admits a genus-one fibration $f:X \to \bP_{k}^1$
by \cref{prop:genus-one-fibration}. 
If $f$ is quasi-elliptic then $X(k)\neq \emptyset$ by \cref{lem:quasiell_implies_ss}, 
and we are done. 
Suppose that $f$ is an elliptic fibration.  
If there exists $t\in \bP^1_k(k)$ such that the fiber $X_t$ is smooth, 
then $X_t(k)\neq \emptyset$ (by, for example, the Weil conjectures), which 
implies that $X(k)\neq \emptyset$, and we are done. 
Suppose that $X_t$ is singular for all $t\in \bP^1_k(k)$. 
Let $\overline{f}:\overline{X}\to \bP_{\overline{k}}^1$ denote the base extension of 
$f$ to $\overline{k}$. 
Note that the fiber $\overline{X}_t$ of $\overline{f}$ over $t\in \bP^1_k(k)$ 
is the same as the base extension of $X_t$ to $\overline{k}$. 
Let $\Sigma$ be the set of closed points $t$ of $\bP_{\overline{k}}^1$ 
such that the fiber $\overline{X}_t$ is singular. 
We have $\bP^1_k(k)\subset \Sigma$.  
By \cite[Corollary 4.1.12]{Cossec-Dolgachev-Liedtke_EnriquesI}, it follows that 
\begin{equation}\label{eq:forEulerChar}
24 \geq \sum_{t\in \Sigma} e(\overline{X}_t) 
\geq \sum_{t\in \bP^1_k(k)} e(\overline{X}_t)
\end{equation}
where $e(\overline{X}_t) 
= \sum_{i\geq 0} (-1)^i \dim_{\bQ_\ell}H^i_{\textrm{\'et}}(\overline{X}_t, \bQ_\ell)$. 
We have $e(\overline{X}_t) < 2$ for some $t\in \bP^1_k(k)$; 
because otherwise \eqref{eq:forEulerChar} would yield 
$24 \geq 2|\bP^1_k(k)| = 2(|k|+1) > 24$ since $|k|\geq 13$.   
Let $t\in \bP^1_k(k)$ be a point such that $e(\overline{X}_t) < 2$. 
Then, by the classification of possible singular fibers 
(e.g. \cite[Theorem 4.1.4]{Cossec-Dolgachev-Liedtke_EnriquesI}), 
it follows that $\overline{X}_t$ is irreducible. 
Since no fiber is multiple (e.g. \cite[Proposition 11-1.6]{Huybrechts_K3book}), 
the singular fiber $X_t$ is a geometrically integral curve of arithmetic genus $1$. 
Hence, its normalization $(X_t)^\nu$ is $\bP^1_k$, and in particular, it has 
a $k$-point, say $r$. Then $\Spec(k)\xrightarrow{r} (X_t)^\nu \to X_t \to X$ is a $k$-point 
of $X$. 
This completes the proof. 
\end{proof}

\section{Examples of pointless K3 surfaces}\label{sec:examples}
This section presents explicit projective models of pointless K3 surfaces over
finite fields. These examples prove \cref{th:pointlessK3surfaces}.
Some of the examples presented below were found with the assistance of generative AI. We independently verified all examples obtained in this way using OSCAR: we checked that they are smooth and pointless.

\subsection{Case \texorpdfstring{$q\leq 5$}{q <= 5}}
\paragraph{Case $q = 2$}
According to the data cited in \cite{Kedlaya-Sutherland2016}, 
up to the action of $\mathrm{PGL}_4(\bF_2)$, 
there exist exactly $34$ quartic K3 surfaces in $\bP^3$ over $\bF_2$ with no $\bF_2$-point.

\paragraph{Case $q = 3$}
Using a brute-force computer search, we found a pointless quartic K3 surface in $\bP^3$ over $\bF_3$, which is defined by the following equation:
\[ \begin{split}
&x_0^4 + x_0^3x_1 + x_0^3x_2 + x_0^3x_3 + x_0^2x_1^2 + x_0^2x_1x_2 + x_0^2x_1x_3 \\
&+ x_0^2x_2^2 +  x_0^2x_2x_3 + x_0^2x_3^2 + x_0x_1^3 + x_0x_1^2x_2 + x_0x_1^2x_3 \\
&+ x_0x_1x_2^2 + x_0x_1x_2x_3 + x_0x_1x_3^2 + x_0x_2^3 + x_0x_2^2x_3 + x_0x_2x_3^2 \\
&+ x_0x_3^3 + x_1^4 + x_1^3x_2 + x_1^3x_3 + x_1^2x_2^2 + 2x_1^2x_2x_3 + x_1^2x_3^2 \\
&+ x_1x_2^2x_3 + 2x_1x_2x_3^2 + x_2^4 + x_2^3x_3 + x_2^2x_3^2 + x_3^4 = 0.
\end{split}
\]

\paragraph{Case $q = 4$}
A computer search assisted by generative AI produced the following example over $\bF_4$.
Note that $\bF_4 = \bF_2[t]/(t^2 + t + 1)$ since the polynomial $t^2 + t + 1$ is irreducible over $\bF_2$. Let $\alpha\in \bF_4$ denote the class of $t$. 
We consider the weighted projective space $\bP_{\bF_4}(3,1,1,1)$ with weighted homogeneous coordinates $(w:x:y:z)$. 
Let $X\subset \bP_{\bF_4}(3,1,1,1)$ be the surface defined by 
the weighted homogeneous polynomial 
\[F(w,x,y,z) \coloneq w^2+A(x,y,z)w+B(x,y,z)
\]
where 
\begin{align*}
A(x,y,z) = z^3 + \alpha yz^2+y^2z + (\alpha+1) y^3 + xz^2 + (\alpha+1) xyz + x^3
\end{align*}
and 
\begin{align*}
B(x,y,z)={}&
 \alpha z^6+yz^5+\alpha y^2z^4+(\alpha+1)y^4z^2 \\
&+(\alpha+1)y^5z+(\alpha+1)y^6+xz^5
  +\alpha xy^3z^2 \\
&+\alpha xy^4z+xy^5+x^2z^4+x^2yz^3 \\
&+\alpha x^2y^2z^2+(\alpha+1)x^2y^3z
  +(\alpha+1)x^2y^4 \\
&+x^3yz^2+x^3y^3+\alpha x^4yz \\
&+\alpha x^5z+x^5y+(\alpha+1)x^6.
\end{align*}
A direct computation shows that $X(\bF_4)=\emptyset$.
Furthermore, the Jacobian criterion shows that $X$ is smooth and hence a K3 surface. 

\paragraph{Case $q = 5$}
The Fermat quartic surface $x_0^4 + x_1^4 + x_2^4 + x_3^4 = 0$ in $\bP^3$
is a pointless K3 surface over $\bF_5$.

\subsection{Case \texorpdfstring{$q\in \{7,8,9,11\}$}{q in {7,8,9,11}}}
The proof of \cref{th:criteria-for-having-pts} (ii) shows that  
if a K3 surface $X$ over $\bF_q$ admitting an elliptic fibration is pointless, then the configuration of its singular fibers is subject to certain constraints. 
Taking these constraints into account, we used generative AI to search for such examples 
and found one over $\bF_q$ for each $q\in \{7,8,9,11\}$.
Each example is given as a double cover $X \to \bP^1\times \bP^1$, and the composition  
$X\to \bP^1\times \bP^1 \to \bP^1$ defines an elliptic fibration. 

Before presenting the examples, we recall the construction. Let $k$ be a field, and set $Y \coloneq \bP^1\times \bP^1$, with homogeneous coordinates $[t_0:t_1]$ and $[x_0:x_1]$ on the two factors. 
Let $F$ and $G$ be bi-homogeneous polynomials in $k[t_0, t_1, x_0, x_1]$ of bi-degrees
$(2,2)$ and $(4,4)$, respectively. 
Let $X$ be the double cover of $Y$ defined by the equation 
\begin{equation}\label{eq:double_cover_of_P1timesP1}
w^2 - F(t_0, t_1, x_0,x_1) w - G(t_0, t_1, x_0,x_1) = 0. 
\end{equation}
Here, $w$ is regarded as a fiber coordinate on the total space of $\cO_Y(2,2)$.
One sees that $X$ is a K3 surface if it is smooth. 
When $\mathrm{char}(k)\neq 2$ and $F = 0$, $X$ is a double cover of $Y$ branched along the curve defined by $G = 0$, and it is the standard construction of a cyclic cover. 
In characteristic $2$, however, 
$F$ needs to be nonzero to obtain a pointless K3 surface over $\bF_8$. 

For an $(n+1) \times (n+1)$ matrix $C = (c_{ij})_{1\leq i,j\leq n+1}$, 
we define a bi-homogeneous polynomial $F_C\in k[t_0, t_1, x_0,x_1]$ 
of bi-degree $(n, n)$ by 
\[ F_C \coloneq \sum_{0\leq i,j\leq n} 
c_{i+1\,j+1} t_0^{n-i} t_1^{i} x_0^{n-j} x_1^{j}.
\]

\paragraph{Case $q=7$}
Set $F = 0$ and $G = F_C$ in \eqref{eq:double_cover_of_P1timesP1}, where 
\[ C=
\begin{pmatrix}
6 & 0 & 3 & 0 & 3\\
2 & 6 & 6 & 4 & 0\\
1 & 2 & 1 & 6 & 6\\
5 & 1 & 1 & 3 & 0\\
6 & 0 & 3 & 0 & 3
\end{pmatrix}. 
\]
Then, the resulting surface $X$ over $\bF_7$ is a pointless K3 surface. 

\paragraph{Case $q=8$}
The field $\bF_8$ can be written as $\bF_8 = \bF_2[\alpha]/(\alpha^3 + \alpha + 1)$. 
Set $F = F_A$ and $G = F_B$ in \eqref{eq:double_cover_of_P1timesP1}, where 
\[ A=
\begin{pmatrix}
1 & 1 & 1\\
1 & \alpha^2 & 1\\
1+\alpha & 1 & 1
\end{pmatrix},
\qquad
B=
\begin{pmatrix}
1+\alpha & 0 & 1 & \alpha & 1+\alpha^2\\
0 & 0 & 0 & \alpha & \alpha^2\\
1 & 0 & \alpha+\alpha^2 & 1+\alpha & 1+\alpha\\
\alpha & \alpha & 0 & \alpha+\alpha^2 & 0\\
\alpha+\alpha^2 & 1+\alpha+\alpha^2 & 1+\alpha^2
    & \alpha^2 & 1+\alpha+\alpha^2
\end{pmatrix}.
\]
Then, the resulting surface $X$ over $\bF_8$ is a pointless K3 surface.

\paragraph{Case $q=9$}
The field $\bF_9$ can be written as $\bF_9 = \bF_3[i]/(i^2 + 1)$. 
Set $F = 0$ and $G = F_C$ in \eqref{eq:double_cover_of_P1timesP1}, where 
\[ C=
\begin{pmatrix}
2+i & 1+i & i & 1+2i & 1+2i\\
1 & 0 & 2i & 1+2i & 0\\
1 & 0 & 1+2i & 0 & i\\
2+i & 0 & 1+i & 2i & 0\\
2+2i & 2+i & 2 & 1+i & 1+i
\end{pmatrix}.
\]
Then, the resulting surface $X$ over $\bF_9$ is a pointless K3 surface. 

\paragraph{Case $q=11$}
Set $F = 0$ and $G = F_C$ in \eqref{eq:double_cover_of_P1timesP1}, where 
\[ C=
\begin{pmatrix}
6 & 8 & 8 & 6 & 2\\
2 & 0 & 2 & 7 & 0\\
1 & 5 & 5 & 1 & 4\\
9 & 0 & 9 & 4 & 0\\
6 & 8 & 8 & 6 & 2
\end{pmatrix}.
\]
Then, the resulting surface $X$ over $\bF_{11}$ is a pointless K3 surface. 

\begin{proof}[Proof of \cref{th:pointlessK3surfaces}]  
For every prime power $q\leq 11$, the examples above provide an explicit 
pointless K3 surface over $\bF_q$.
\end{proof}

\bibliographystyle{alpha}
\bibliography{literature}

\begin{thebibliography}{MRW08}

\bibitem[AEKS]{AEKS}
Asher Auel, Philip Engel, Kiran Kedlaya, and Sair Shaikh.
\newblock Pointless {K}3 surfaces over finite fields.
\newblock Preprint.

\bibitem[BC22]{brandhorst-cattaneo}
Simon Brandhorst and Alberto Cattaneo.
\newblock {Prime Order Isometries of Unimodular Lattices and Automorphisms of IHS Manifolds}.
\newblock {\em International Mathematics Research Notices}, 10 2022.
\newblock rnac279.

\bibitem[BF25]{Bernasconi-Stefano2025}
Fabio Bernasconi and Stefano Filipazzi.
\newblock Rational points on 3-folds with nef anti-canonical class over finite fields.
\newblock {\em Proc. Lond. Math. Soc. (3)}, 130(1):Paper No. e70014, 31, 2025.

\bibitem[BH23]{brandhorst-hofmann}
Simon Brandhorst and Tommy Hofmann.
\newblock Finite subgroups of automorphisms of {{\(K3\)}} surfaces.
\newblock {\em Forum Math. Sigma}, 11:57, 2023.
\newblock Id/No e54.

\bibitem[BH26]{Brandhorst-Hofmann2026}
Simon Brandhorst and Tommy Hofmann.
\newblock Automorphisms and canonical forms of integer lattices, 2026.

\bibitem[BLvL20]{Brigh-Logan-vanLuijk2020}
Martin Bright, Adam Logan, and Ronald van Luijk.
\newblock Finiteness results for {K}3 surfaces over arbitrary fields.
\newblock {\em Eur. J. Math.}, 6(2):336--366, 2020.

\bibitem[BM]{Muller-Brandhorst2026:db}
Simon Brandhorst and Stevell Muller.
\newblock Deflatdb a database of genera of definite integer lattices.

\bibitem[BT26]{ancillary-files}
Simon Brandhorst and Yuta Takada.
\newblock No pointless {K}3 surfaces for {$q>11$}.
\newblock [Dataset] \href{https://doi.org/10.5281/zenodo.23696784}{https://doi.org/10.5281/zenodo.23696784}, October 2026.

\bibitem[CDL25]{Cossec-Dolgachev-Liedtke_EnriquesI}
Fran\c{c}ois Cossec, Igor Dolgachev, and Christian Liedtke.
\newblock {\em Enriques surfaces. {I}}.
\newblock Springer, Singapore, second edition, [2025] \copyright 2025.
\newblock With an appendix by S. Kond\=o.

\bibitem[Coh15]{cohn-bounds}
Henry Cohn.
\newblock Table of sphere packing density bounds.
\newblock \url{https://hdl.handle.net/1721.1/153311}, 2024-01-15.

\bibitem[FHS19]{Fieker-Hofmann-Sircana2019}
Claus Fieker, Tommy Hofmann, and Carlo Sircana.
\newblock On the construction of class fields.
\newblock In {\em ANTS XIII. Proceedings of the thirteenth algorithmic number theory symposium, University of Wisconsin-Madison, WI, USA, July 16--20, 2018}, pages 239--255. Berkeley, CA: Mathematical Sciences Publishers (MSP), 2019.

\bibitem[Gou20]{Gouvea_padic}
Fernando~Q. Gouv\^ea.
\newblock {\em {$p$}-adic numbers}.
\newblock Universitext. Springer, Cham, third edition, 2020.
\newblock An introduction.

\bibitem[Huy16]{Huybrechts_K3book}
Daniel Huybrechts.
\newblock {\em Lectures on {K}3 surfaces}, volume 158 of {\em Cambridge Studies in Advanced Mathematics}.
\newblock Cambridge University Press, Cambridge, 2016.

\bibitem[Jac62]{jacobowitz}
Ronald Jacobowitz.
\newblock Hermitian forms over local fields.
\newblock {\em Amer. J. Math.}, 84:441--465, 1962.

\bibitem[Kaw24]{Kawabe2024}
Daiki Kawabe.
\newblock Chow motives of quasi-elliptic surfaces.
\newblock {\em Osaka J. Math.}, 61(4):559--567, 2024.

\bibitem[Kir16]{KirschmerHabil}
Markus Kirschmer.
\newblock {\em Definite quadratic and hermitian forms with small class number}.
\newblock Habilitationsschrift, RWTH Aachen University, Aachen, Germany, 2016.
\newblock \href{http://www.math.rwth-aachen.de/~Markus.Kirschmer/papers/herm.pdf}{http://www.math.rwth-aachen.de/~Markus.Kirschmer/papers/herm.pdf}.

\bibitem[Kir19]{Kirschmer2019}
Markus Kirschmer.
\newblock Determinant groups of {H}ermitian lattices over local fields.
\newblock {\em Arch. Math. (Basel)}, 113(4):337--347, 2019.

\bibitem[Kon20]{Kondo_K3book}
Shigeyuki Kond\=o.
\newblock {\em {$K3$} surfaces}, volume~32 of {\em EMS Tracts in Mathematics}.
\newblock EMS Publishing House, Berlin, [2020] \copyright 2020.
\newblock Translated from the Japanese original by the author.

\bibitem[KS16]{Kedlaya-Sutherland2016}
Kiran~S. Kedlaya and Andrew~V. Sutherland.
\newblock A census of zeta functions of quartic {K}3 surfaces over {$\Bbb F_2$}.
\newblock {\em LMS J. Comput. Math.}, 19:1--11, 2016.

\bibitem[LM18]{Liedtke-Matsumoto2018}
Christian Liedtke and Yuya Matsumoto.
\newblock Good reduction of {K}3 surfaces.
\newblock {\em Compos. Math.}, 154(1):1--35, 2018.

\bibitem[McM11]{McMullen2011}
Curtis~T. McMullen.
\newblock K3 surfaces, entropy and glue.
\newblock {\em J. Reine Angew. Math.}, 658:1--25, 2011.

\bibitem[Mil69]{Milnor1969}
John Milnor.
\newblock On isometries of inner product spaces.
\newblock {\em Invent. Math.}, 8:83--97, 1969.

\bibitem[MRW08]{mossinghoff:minimal-mahler}
Michael~J. Mossinghoff, Georges Rhin, and Qiang Wu.
\newblock Minimal {Mahler} measures.
\newblock {\em Exp. Math.}, 17(4):451--458, 2008.

\bibitem[{Nik}80]{nikulin1980_forms}
Viacheslav~V. {Nikulin}.
\newblock {Integral symmetric bilinear forms and some of their applications}.
\newblock {\em {Math. USSR, Izv.}}, 14:103--167, 1980.

\bibitem[Nyg83]{Nygaard1983}
N.~O. Nygaard.
\newblock The {T}ate conjecture for ordinary {$K3$}\ surfaces over finite fields.
\newblock {\em Invent. Math.}, 74(2):213--237, 1983.

\bibitem[OSC25]{OSCAR}
O{SCAR} -- {O}pen {S}ource {C}omputer {A}lgebra {R}esearch system, {V}ersion 1.6.0, 2025.

\bibitem[Poo17]{Poonen_rational-points}
Bjorn Poonen.
\newblock {\em Rational points on varieties}, volume 186 of {\em Graduate Studies in Mathematics}.
\newblock American Mathematical Society, Providence, RI, 2017.

\bibitem[PS97]{plesken-souvignier}
Wilhelm {Plesken} and Bernd {Souvignier}.
\newblock {Computing isometries of lattices}.
\newblock {\em {J. Symb. Comput.}}, 24(3-4):327--334, 1997.

\bibitem[PZ25]{Patakfalvi-Zdanowicz2019}
Zsolt Patakfalvi and Maciej Zdanowicz.
\newblock The weak {B}eauville--{B}ogomolov decomposition in characteristic $p\geq 0$, 2025.
\newblock arXiv:1912.12742v3.

\bibitem[{Sta}26]{StacksProject}
The {Stacks project authors}.
\newblock The stacks project.
\newblock \url{https://stacks.math.columbia.edu}, 2026.

\bibitem[Tae20]{Taelman20}
Lenny Taelman.
\newblock Ordinary {K}3 surfaces over a finite field.
\newblock {\em J. Reine Angew. Math.}, 761:141--161, 2020.

\bibitem[vzGG13]{modern-computer-algebra}
Joachim von~zur Gathen and J{\"u}rgen Gerhard.
\newblock {\em Modern computer algebra}.
\newblock Cambridge: Cambridge University Press, 3rd edition, 2013.

\bibitem[Yu12]{Yu2012}
Jeng-Daw Yu.
\newblock Special lifts of ordinary {K}3 surfaces and applications.
\newblock {\em Pure Appl. Math. Q.}, 8(3):805--824, 2012.

\end{thebibliography}
\end{document}